\documentclass[11pt]{article}
\usepackage{amsmath,amsthm}
\usepackage{amsfonts}
\usepackage{amssymb}
\usepackage{graphicx}
\usepackage{authblk}
\usepackage{color}
\usepackage{epstopdf}
\usepackage[symbol]{footmisc}

\newlength{\spacing}
\newtheorem{lemma}{Lemma}[section]

\newtheorem{remark}{Remark}[section]

\numberwithin{equation}{section}

\def\beq{\begin{equation}}
\def\eeq{\end{equation}}

\title{A generalized framework to construct third order WENO weights using weight limiter functions}
\author[1]{Ritesh Kumar Dubey }
\author[2]{Sabana Parvin }
\affil[*]{Research Institute \& Department of Mathematics, SRM Institute of Science and Technology, Chennai, India}
\affil[1]{\textit { riteshkumar.d@res.srmuniv.ac.in}}
\affil[2]{\textit {sabanaparvin.s@res.srmuniv.ac.in}}

\date{}
\begin{document}

	\maketitle
	
\begin{abstract}
The main aim of this work is not to improve any existing non-linear weight but to give a generalized framework for the construction of non-linear weights to get non-oscillatory third order WENO schemes. It is done by imposing necessary conditions on weights to get non-oscillatory WENO reconstruction which give further insight on the  structure of weights to ensure non-occurrence of oscillations and characterize the solution region for third order accuracy. This framework for WENO weights is new and completely different from the prevailing existing approach. New non-linear weights are designed using a function of smoothness parameter termed as weight limiter functions. Many such weight limiter functions are given and analyzed. These new weights are simple and by construction guarantee for exact third order accuracy in smooth solution region including smooth extrema away from critical point. Numerical results for various test problems are given and compared. Results show that proposed weights give third order accuracy without loosing the non-oscillatory shock capturing ability of the resulting scheme.
\end{abstract} 

{\bf Keywords:}	Hyperbolic  conservation  laws, Third order WENO reconstructions, Limiter functions, Non-linear weights, Data dependent stability.\\
{\bf AMS subject classifications}. 65M06, 65M06, 35L65
	
	 \section{ Introduction}
    The initial value problem for system of conservation laws which models flow phenomena in Gas Dynamics, Aero-dynamics, Astrophysical modeling, meteorology and weather prediction etc. can be written in one space dimension with given $f(u):\mathbb R^m\to \mathbb R^m~ \text{and}~u_0:\mathbb R\to \mathbb R^m$ as
    \begin{equation}\label{eqn1}
    \begin{aligned} 
    u_t+f(u)_x&=0,~~~~~~~x\in \mathbb R\times(0,\infty)\\ 
    u&=u_0 (x),~~x\in \mathbb R\times(t=0)
    \end{aligned}
    \end{equation}
   where $u:\mathbb R\times(0,\infty)\to \mathbb R^m$ is the unknown conserved variable.
   In general the closed form solution is not known for such complex systems which normally have strongly irregular solutions (large jumps, discontinuities) and such problems may involve complicated smooth solution region structures. The numerical treatment of such problems faces typical complexity, for example, oscillations near shock and contact discontinuity which makes a high order numerical scheme unstable. On the other hand in certain situations, the time of evolution of these complex structures are so long that it is impractical to use low order numerical methods to obtain an acceptable resolution. Therefore, development and analysis of high order accurate and non-oscillatory shock capturing numerical schemes for problem \eqref{eqn1} has been an active area of research. Among existing non-oscillatory shock capturing schemes, total variation diminishing (TVD) and weighted essentially non-oscillatory schemes (WENO) have been of great interest because of their ability to capture the discontinuity with high resolution. 
   \par The concept of TVD schemes is introduced first in \cite{harten1983} and was further utilized to develop and analyze high resolution TVD schemes in \cite{Sweby1984, MUNZ198818,toro2000}. The basic idea therein is to use flux limiters to get high order reconstructed value of the flux. More details on flux limited TVD methods and applications can be found in \cite{ZHANG2015114, DUBEY2013325}. Unfortunately, though TVD schemes ensure for removal of spurious oscillations from the numerical approximations however they are criticized due to their degenerate accuracy at extrema \cite{tadmor1988}. This degeneracy causes cornered approximation to smooth solution due to clipping error.

\par On the other hand, essentially non-oscillatory (ENO) scheme based on cell-average is proposed by Harten et al. in \cite{harten1987uniformly}. The idea behind ENO schemes is to choose the smoothest stencil among several candidate stencil for high order accurate approximation of the flux function at cell boundaries and also to avoids oscillations near shocks. However,  the cell-average approach to reconstruct point values from given cell average values in multidimensional case is computationally costly compared to the flux version of efficient ENO scheme (ENO schemes based on point values) given in \cite{shu1988efficient,shu1989efficient}. Later many improved version of ENO schemes has been developed see \cite{shu1990numerical}. The most popular of them is cell averaged weighted ENO scheme (WENO  introduced by Liu, Osher and Chan \cite{liu1994weighted}. It uses non-linear convex combination of all interpolating polynomial obtained from candidate stencils of ENO scheme, which results in to a higher order accurate non-oscillatory scheme compared to ENO scheme using the same stencil. The most significant contribution of this technique is the construction of non-linear weights and smoothness indicator based on undivided differences. Later Jiang and Shu \cite{jiang1996efficient} introduced finite difference flux version WENO scheme popularly named as WENO-JS scheme by modifying the smoothness measurement and extended the scheme up to $5th$ order accuracy. Later, Henrick et al. found that non-linear weights in WENO-JS failed to recover optimal (ideal) order of accuracy at critical points and they presented mapped WENO scheme (WENO-M) \cite{henrick2005mapped}. In mapped WENO the non-linear weights are constructed through the construction of mapping function so that weights remain as close as possible to optimal (Ideal) weights except at highly non-smooth regions. Another approach was developed by Borges et al. \cite{borges2008improved} by using global smoothness measurement of fifth order WENO scheme named as WENO-Z scheme, having the same accuracy of WENO-M scheme but with less computational cost. Further many modified versions of WENO schemes have been developed to improve the order of accuracy or to get computationally cheaper results by changing the smoothness indicator of the non-linear weights \cite{castro2011high,ha2013improved,rathan2018}. In \cite{yamaleev2009third} energy stable WENO3 schemes are proposed.
\par Recently active research work has been carried out specific to improvisation of non-linear weights of third order WENO schemes.  The classical WENO-JS non-linear weights are modified to achieve optimal third order accuracy regardless of critical point in \cite{DON2013347}. These schemes are known as WENO-Z3 schemes. Further, improvement to get faster WENO schemes, WENO-Z3 wights are modified in \cite{xiaoshuai2015high} termed as WENO-N3 which are again modified in \cite{wu2016new} to get third order accuracy at critical points and termed and WENO-NP3. In \cite{rathan2017}, WENO-NP3 weights are modified by changing the smoothness indicator named as WENO-F3. Very recently, in \cite{xu2018improved}, WENO-Z3 weights are significantly modified to get ENO solution near strong discontinuity and termed as WENO-P+3 schemes. It is also shown that WENO-P+3 scheme outperforms WENO3 schemes using existing improved non-linear weights.

\par The main aim of this work is not to improve any existing weight but to present a generic framework to devise non-linear WENO weights which guarantee by construction for third order accuracy in smooth solution region. This is achieved via analyzing the weights of third order WENO reconstruction in the light of data dependent stability notion given in \cite{dubey2016data, dubey2017investigation} along with the flux limiters framework.  More precisely, simple non-linear weights are designed by using weight limiter as a function of a smoothness measurement which is again a function of consecutive gradient ratio. An important feature is that weight limiter function controls the weights such that they achieve ideal weight for third order accuracy in smooth region including most of extreme points and give non-oscillatory approximation for discontinuity. Various examples of such weight limiter functions are given.
   
   The organization of the paper is as follows: Section \ref{sec2} provides brief review of WENO scheme  and in particular third order WENO weights in section \ref{sec2p1}. The main detailed contribution of work is given section \ref{sec3}. In section \ref{sec4} demonstrates the computational testing and quantitative performance of the third order WENO scheme using proposed non-linear weights. Conclusion and remark on ongoing work given in section \ref{sec5}.
     
\section{ Review of WENO scheme}\label{sec2}
In this section we briefly describe the flux version third order WENO finite difference scheme for one dimensional scalar conservation law proposed in \cite{jiang1996efficient}. Let $I_i$ be a partition of a given domain with $i^{th}$ cell $I_i=\left[x_{i-\frac{1}{2}},x_{i+\frac{1}{2}}\right] $, center of $I_i$ is $\displaystyle x_i=\frac{x_{i-\frac{1}{2}}+x_{i+\frac{1}{2}}}{2}$ and function value $f$ at node $x_i$ is given by $f_i=f(x_i)$. For simplicity $\{x_{i+\frac{1}{2}}\}_{i}$ is uniformly spaced, notation $\Delta x=x_{i+\frac{1}{2}}-x_{i-\frac{1}{2}}$ indicate the size of $I_i$. Moreover, notation $u_i^n$ used for approximation to $u$ at the grid point $(x_i,t^n)$ where $t^n=n\Delta t$ is the discrete value in the time direction.\\
The semi-discretized approximation of the one -dimensional hyperbolic conservation law \eqref{eqn1} can be expressed as:
\begin{gather}\label{eqn2.1}
\frac{du_i(t)}{dt}=-\frac{1}{\Delta x}\left(\hat{f}_{i+\frac{1}{2}}-\hat{f}_{i-\frac{1}{2}}\right)=:L(u)
\end{gather}
where $u_i(t)$ is the numerical approximation to the point values $u(x_i,t)$ and the numerical flux $\hat{f}$ is a function of $(r+s)$ variables i.e., $\hat{f}_{i+\frac{1}{2}}=\hat{f}(u_{i-r},...,u_{i+s}) $. The numerical flux function should satisfy Lipschitz continuity in each of its arguments, be consistent with the physical flux $f$, that is, $\hat{f}(u,u,...,u)=f(u)$.\\
To compute numerical flux $\hat{f}_{i\pm\frac{1}{2}}$, a function $h$ is defined implicitly by the following equation (see Lemma $2.1$ of \cite{shu1989efficient}).
\begin{align}\label{eqn2.2}
f(u(x))=\frac{1}{\Delta x}\int_{x-\frac{\Delta x}{2}}^{x+\frac{\Delta x}{2}}h(\xi) d\xi.
\end{align}
differentiation of the above equation \eqref{eqn2.2} at the point $x=x_i$ yields,
\begin{align}\label{eqn2.3}
\frac{\partial f}{\partial x}\mid_{x=x_i}=\frac{1}{\Delta x}\left(h_{i+\frac{1}{2}}-h_{i-\frac{1}{2}}\right)
\end{align}
which indicates that the numerical flux $\hat{f}$ approximates $h$ at cell boundaries $x_{i\pm\frac{1}{2}}$ with high order accuracy, that is, 
\begin{equation}\nonumber
\hat{f}_{i\pm\frac{1}{2}}=h_{i\pm\frac{1}{2}}+O(\Delta x^r)
\end{equation}
with $r$ depending on the degree of interpolation.\\
Using equation \eqref{eqn2.3} in equation \eqref{eqn2.1}, we have
\begin{equation}\label{eqn2.4}
\frac{du_i(t)}{dt}=-\frac{1}{\Delta x}\left(h_{i+\frac{1}{2}}-h_{i-\frac{1}{2}}\right)\approx -\frac{1}{\Delta x}\left(\hat{f}_{i+\frac{1}{2}}-\hat{f}_{i-\frac{1}{2}}\right)
\end{equation}
In order to avoid entropy violating solution  and ensure numerical stability flux $f(u)$ is splitted into the two components $f^{+}$ and $f^{-}$ such that 
\begin{align}\label{eqn2.5}
f(u)=f^{+}(u)+f^{-}(u)
\end{align}
where $\frac{df^{+}(u)}{du}\ge 0$ and $\frac{df^{-}(u)}{du}\le 0$. Among many flux splitting methods, the following global Lax-Friedrichs splitting is heavily used for its simplicity and capability to produce very smooth fluxes, which is defined as
\begin{align}
f^{\pm}(u)=\frac{1}{2}(f(u)\pm \alpha u)
\end{align}
where $\alpha =\max_{u} |f^{\prime}(u)|$. Let $\hat{f}^{+}_{i+\frac{1}{2}}$ and $\hat{f}^{-}_{i+\frac{1}{2}}$ be the numerical fluxes obtained from the positive and negative parts of $f(u)$ respectively and from \eqref{eqn2.5} we have 
\begin{align}\label{eqn2.7}
\hat{f}_{i+\frac{1}{2}}=\hat{f}^{+}_{i+\frac{1}{2}}+\hat{f}^{-}_{i+\frac{1}{2}}
\end{align}
In the following approximation of $\hat{f}^{+}_{i+\frac{1}{2}}$ is describe for third order accuracy. Note that   negative part $\hat{f}^{-}_{i+\frac{1}{2}}$ of \eqref{eqn2.7} can be approximated accordingly as it is symmetric to positive part $\hat{f}^{+}_{i+\frac{1}{2}}$ with respect to $x_{i+\frac{1}{2}}$. Also, for simplicity we will drop the $'+'$  sign in the superscript.

\subsection{ Third order WENO}\label{sec2p1}
The construction of $\hat{f}^{+}_{i+\frac{1}{2}}$	for the classical third-order WENO scheme uses a $3$-point stencil $S(i):=\{x_{i+k-1},x_{i+k},x_{i+k+1}\},(k=0,1)$ which is subdivided into two candidate sub-stencils. Let 
\begin{equation}
S_k(i):=\{x_{i+k-1},x_{i+k}\},~~~k=0,1
\end{equation}
be the sub-stencil consisting of $2$-points starting at $x_{i+k-1}$ and let 
	\begin{equation}
	\hat{f}^{k}_{i+\frac{1}{2}}=\sum_{l=0}^{1} c_{k,l} f_{i+k-1+l} 
	\end{equation}
	be the first degree polynomial approximation constructed on the stencil $S_k(i)$ to approximate the value $h(x_{i+\frac{1}{2}})$ where $c_{k,l},~l=0,1$ are the Lagrange interpolation coefficients depending on the left-shift parameter $k$. The expression of $\hat{f}^{k}_{i+\frac{1}{2}}$ can be written as
	\begin{equation}\label{eqn2.10}
	\begin{aligned}
	\hat{f}^{0}_{i+\frac{1}{2}}=\frac{3}{2}f_i-\frac{1}{2}f_{i-1}\\
	\hat{f}^{1}_{i+\frac{1}{2}}=\frac{1}{2}f_i+\frac{1}{2}f_{i+1}
	\end{aligned}
	\end{equation}
	To define $\hat{f}^{j}_{i-\frac{1}{2}}$ each index needs to be shifted by $-1$. Moreover, the Taylor expansion of the equations in \eqref{eqn2.10} are given by
	\begin{equation}\label{eqn2.11}
	\begin{aligned}
	\hat{f}^{0}_{i+\frac{1}{2}}=h_{i+\frac{1}{2}}-\frac{\Delta x^2}{4} f^{(2)}(0)+O(\Delta x^3)\\
	\hat{f}^{1}_{i+\frac{1}{2}}=h_{i+\frac{1}{2}}+\frac{\Delta x^2}{4} f^{(2)}(0)+O(\Delta x^3)
	\end{aligned}
	\end{equation}
	These functions are combined to define a new WENO approximation to the value $h_{i+\frac{1}{2}}$, that is,
	\begin{equation}\label{eqn2.12}
	\hat{f}_{i+\frac{1}{2}}=\sum_{k=0}^{1}\omega_k \hat{f}^{k}_{i+\frac{1}{2}}
	\end{equation}
	where non-linear weights $\omega_k$ satisfy the following convexity property
	\begin{equation}\label{convex}
	\sum_{k=0}^1 \omega_k=1,\; \omega_k\geq 0. 
	\end{equation}
	To construct the non-linear weights $\omega _k$ it is first considered the case that the function $h$ is is smooth on all the stencil $S_k(i)$ with $k=0,1$. After that the constants $d_k$ are found such that its linear combination with $\hat{f}^{k}_{i+\frac{1}{2}}$ retains the third order convergence order to $h_{i+\frac{1}{2}}$, that is,
	\begin{equation}
	h_{i+\frac{1}{2}}=\sum_{k=0}^{1} d_k \hat{f}^{k}_{i+\frac{1}{2}}+O(\Delta x^3)
	\end{equation}
	The coefficients $d_k$ are called the ideal or linear weights.The specific values of $d_k$ are as follows \cite{shu1998essentially}
	\begin{equation}
	d_0=\frac{1}{3},~~d_1=\frac{2}{3}
	\end{equation}
	Note that each $d_k\ge 0$ and $\sum_{k=0}^{1} d_k=1$. The non-linear weights $\omega _k$ in \eqref{eqn2.12} are constructed such that final reconstruction become essentially non-oscillatory using following designing criteria 
\begin{itemize}
\item[i.] In smooth regions the non-linear weights should converge to the linear weights with required order of accuracy.
\item[ii.] Weight corresponding to non-smooth or discontinuous region should tend to zero so that contribution from the non-smooth regions in the approximation of $\hat{f}_{i+\frac{1}{2}}$ is negligible
\end{itemize}

\subsubsection{WENO-JS3 Weights \cite{jiang1996efficient}} 
The non-linear WENO-JS3 weights are defined as 
\begin{equation}\label{w-js}
\omega_k=\frac{\alpha_k}{\sum_{l=0}^{1} \alpha_l},~ \alpha_k=\frac{d_k}{(\epsilon+\beta_k)^p}
\end{equation}
where $\epsilon$ is a positive small number which is set to be $\epsilon=10^{-6}$ to avoid division by zero, $p=2$ is chosen to increase the difference of scales of distinct weights ar non-smooth parts of the solution. Note that $\alpha_k$ are the unnormalized weights and $\omega_k$ are the normalized weights. In case of third order accuracy the smoothness indicator  $\beta_k$ can be defined as 
\begin{equation}\label{eqn2.22}
\beta_k=\Delta x \int_{x_{i-\frac{1}{2}}}^{x_{i+\frac{1}{2}}}\left(\frac{d\hat{f}^{k}}{dx}\right)^2 dx, \; k=0,1.
\end{equation}
which reduces to 
\begin{equation}\label{bk}
\begin{aligned}
\beta_0=(f_i-f_{i-1})^2\\
\beta_1=(f_{i+1}-f{i})^2
\end{aligned}
\end{equation}
\subsubsection{WENO-Z3 Weights \cite{don2013accuracy}}
The WENO-JS3 weights \eqref{w-js} were further modified to WENO-Z3. These WENO-Z3 weights are obtained by modifying $\alpha_k$ in \eqref{w-js} as
\begin{equation}\label{yc-weight}
\alpha_k= d_k(1+\frac{\tau}{(\epsilon + \beta_k)})
\end{equation}
where $\tau=|\beta_0-\beta_1|$  and $\beta_k$ are given by \eqref{bk}. 
\subsubsection{WENO-N3 Weights \cite{xiaoshuai2015high}}
A more high resolution and efficient WENO-N3 weights compared to WENO-Z3 are proposed, which are obtained by using
\begin{equation}
\tau=\left|\frac{\beta_0+\beta_1}{2}-\beta_3\right|,
\end{equation}
 where $\beta_3=\frac{13}{12}(f_{i-1}-2f_i+f_{i+1})^2+\frac{1}{4}(f_{i-1}-f_{i+1})^2$.
 WENO-N3 weights are further improved to WENO-NP3 and WENO-F3 weights to achieve third order accuracy at critical points in \cite{wu2016new,rathan2017}. Though WENO3 scheme using these weights give improved approximation to smooth extrema they exhibits oscillations near discontinuity (See Results of WENO-NP3 in figure \ref{FigCompare}b)).  
\subsubsection{WENP-P+3 Weights \cite{xu2018improved}} In 2018, to improve further resolution of WENO3 scheme near strong  discontinuity  WENO-N3 weights are modified by defining 
\begin{equation}\label{WENOp3}
\alpha_k= d_k\left(1+\frac{\tau_p}{(\epsilon + \beta_k)}+\lambda \frac{\beta_k + \epsilon}{\tau_p+\epsilon}\right)
\end{equation}
where $\lambda=\Delta_x^{\frac{1}{6}}$ and $\tau_p=\left|\frac{\beta_0+\beta_1}{2}-\frac{1}{4}(f_{i-1}-f_{i+1})^2\right|.$ Note that though WENO-P+3 scheme outperformed other WENO3 schemes in capturing the strong discontinuity but degenerate to first order accuracy at critical points \cite{xu2018improved}. Also the shock capturing non-oscillatory behavior of WENP-P+3 scheme heavily depends on choice of parameter $\lambda$ (see Figure 5 in \cite{xu2018improved}) and parameter $\epsilon$ in \eqref{WENOp3} see Figure \ref{FigCompare}(a).  


\section{ Construction of new non-linear weights }\label{sec3}
In this section a new approach is given to construct non-linear weights such that the weighted reconstruction \eqref{eqn2.12} achieve essentially non-oscillatory property.  Let us define the parameter of the ratio of consecutive gradients as
\begin{equation}\label{smoothnessparameter}
r_{i}= \frac{\Delta_{-}f_i}{\Delta_{+}f_i}.
\end{equation}
It can be observe from \eqref{eqn2.10} that in case of linear flux i.e., $f=a u,\; a>0$ the reconstructed value at cell interface is a convex combination of second order upwind and centered flux. More precisely in this case \eqref{eqn2.12} reduces to
\begin{equation}\label{reconst}
\hat u_{i+\frac{1}{2}}= \omega_0 \hat u^{0}_{i+\frac{1}{2}} + \omega_1 \hat u^{1}_{i+\frac{1}{2}}
\end{equation}
where $\omega_0\geq 0, \omega_1 \geq 0\; \text{and}\;  \omega_0+\omega_1=1$ . The upwind and centered flux are  
\begin{equation}\label{up2}
\hat u^{0}_{i+\frac{1}{2}}=\frac{3}{2}u_i -\frac{1}{2}u_{i-1}.
\end{equation}
\begin{equation}\label{cent2}
 \hat u^{1}_{i+\frac{1}{2}}= \frac{1}{2}u_i +\frac{1}{2}u_{i+1}.
 \end{equation} 
In order to characterize the weights $\omega$ such that the reconstruction \eqref{reconst} be non-oscillatory, we recall our own results on data dependent stability given in \cite{dubey2016data,dubey2017investigation}
\begin{lemma}\label{lem2}
For the linear transport problem \eqref{eqn4.1} the Forward in time scheme using second order centered flux \eqref{cent2} is data dependent stable and non-oscillatory in the solution data region where $ r_i \in \mathcal{R}_{c}=(-\infty, -1)\cup \left[\frac{a\lambda}{2-a\lambda}, \infty\right)$ for every $i$  where $a\lambda =a\frac{\Delta t}{\Delta x}\leq 1.$ 	
\end{lemma}
\begin{lemma}\label{lem3}
	For the linear transport problem \eqref{eqn4.1}, the forward in time scheme using second order upwind flux \eqref{up2} is data dependent stable and non-oscillatory in the solution data region where $ r_{i-1} \in \mathcal{R}_{up} = \left[-\frac{ 2-3a\lambda}{a\lambda}, 3\right)$ for every $i$ where $a\lambda =a\frac{\Delta t}{\Delta x}\leq \frac{1}{2}.$ 	
\end{lemma}
On dropping out the subscript $i$, following can be deduced.
\begin{itemize}
	\item []{\bf Characterization of centred weight $\omega_1$:} Lemma \ref{lem2} gives the characteristics to be satisfied by the weight $\omega_1$ associated with the centered flux \eqref{cent2} . More precisely, note that for $ a\lambda \rightarrow 0 $ the non-oscillatory region $\mathcal{R}_c$ of centred flux  $(-\infty, -1)\cup \left[\frac{a\lambda}{2-a\lambda}, \infty\right) \rightarrow (-\infty,-1)\cup(0,\infty)$. It clearly concludes that centered flux is oscillatory at $r=0$ irrespective of choice of $\lambda$. Thus in order to avoid contribution of centered flux to construct non oscillatory third order WENO scheme using \eqref{reconst}, a necessary condition is $\omega_1 \rightarrow 0 \; \text{for}\; r \rightarrow 0$. However for $r$ away from zero, $\omega_1$ can be nonzero i.e,  $\omega_1\rightarrow 1 \; \text{for}\; \rightarrow r \pm \infty.$ In particular, under CFL number $a\lambda= \frac{1}{2}$\footnote[1]{Linear stability condition for second order scheme using upwind flux \eqref{up2}}, since  $\mathcal{R}_c =(\infty,-1)\cup(\frac{1}{3},\infty)$, therefore centered weight must be defined such that $\omega_1\rightarrow 0$ as $r\rightarrow 0$ in the interval $[0,1/3)$ however $\omega_1$ can take any value in $[0,1]$ for $\displaystyle r\geq \frac{1}{3}$. \\
\item[]	{\bf Characterization of upwind weight $\omega_0$:}  On the other hand, Lemma \ref{lem3} charaterize the weight $\omega_0$ associated with the upwind flux \eqref{up2}. Note that for $\lambda\rightarrow 0$ non-oscillatory stability region of upwind flux \eqref{up2} i.e., $\mathcal{R}_{up}=   \left[-\frac{ 2-3a\lambda}{a\lambda}, 3\right)\rightarrow (-\infty, 3)$ which suggests that weight $\omega_0\rightarrow 0$ for $r\rightarrow \infty$. In particular, under CFL number $a\lambda= \frac{1}{2}$, since $\mathcal{R}_{up} =[-1,3)$ upwind weight $\omega_0$ can take any value in $[0,1]$ for $r\in [-1,3]$.  
 
\end{itemize}  
\par Thus, necessary conditions for non-oscillatory approximation are $\omega_1\rightarrow 0$ for $r\rightarrow 0$ and $\omega_0\rightarrow 0$ for $r$ away from $[-1,3]$. In particular, for smooth region of solution (which corresponds to $r\approx 1$), convexity property \eqref{convex} along with Lemma \ref{lem2} and \ref{lem3} enables to choose ideal weights $\omega_0=\frac{1}{3}, \; \omega_1=\frac{2}{3}$  for third order accuracy \cite{shu1998essentially}. Moreover, it also follows that third order non-oscillatory approximation can be achieved for $1/3<r<3$ which is the common stability region of centred and upwind flux. 
\par
\begin{figure}[h]
	\centering
	\includegraphics[scale=0.5]{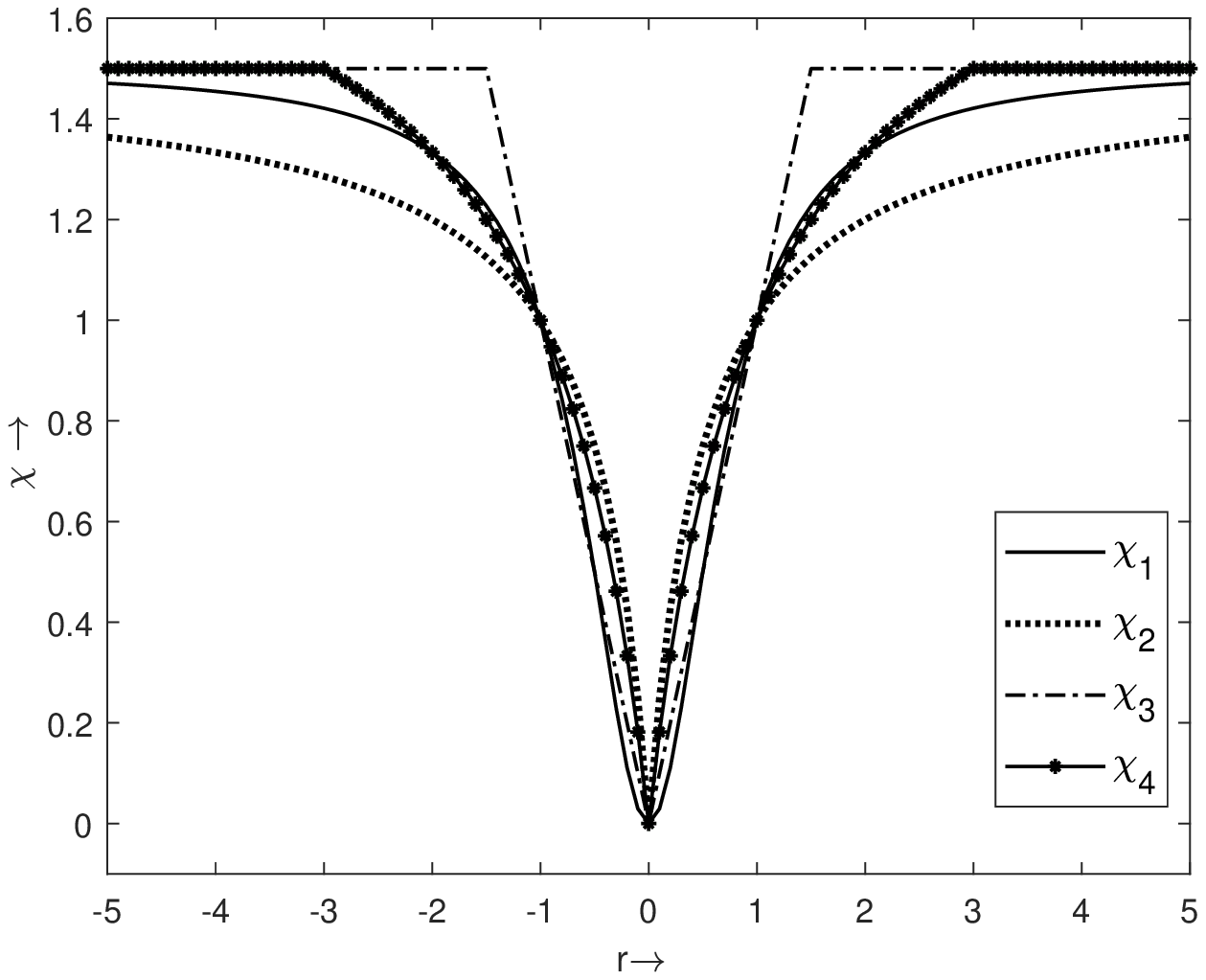}~\includegraphics[scale=0.5]{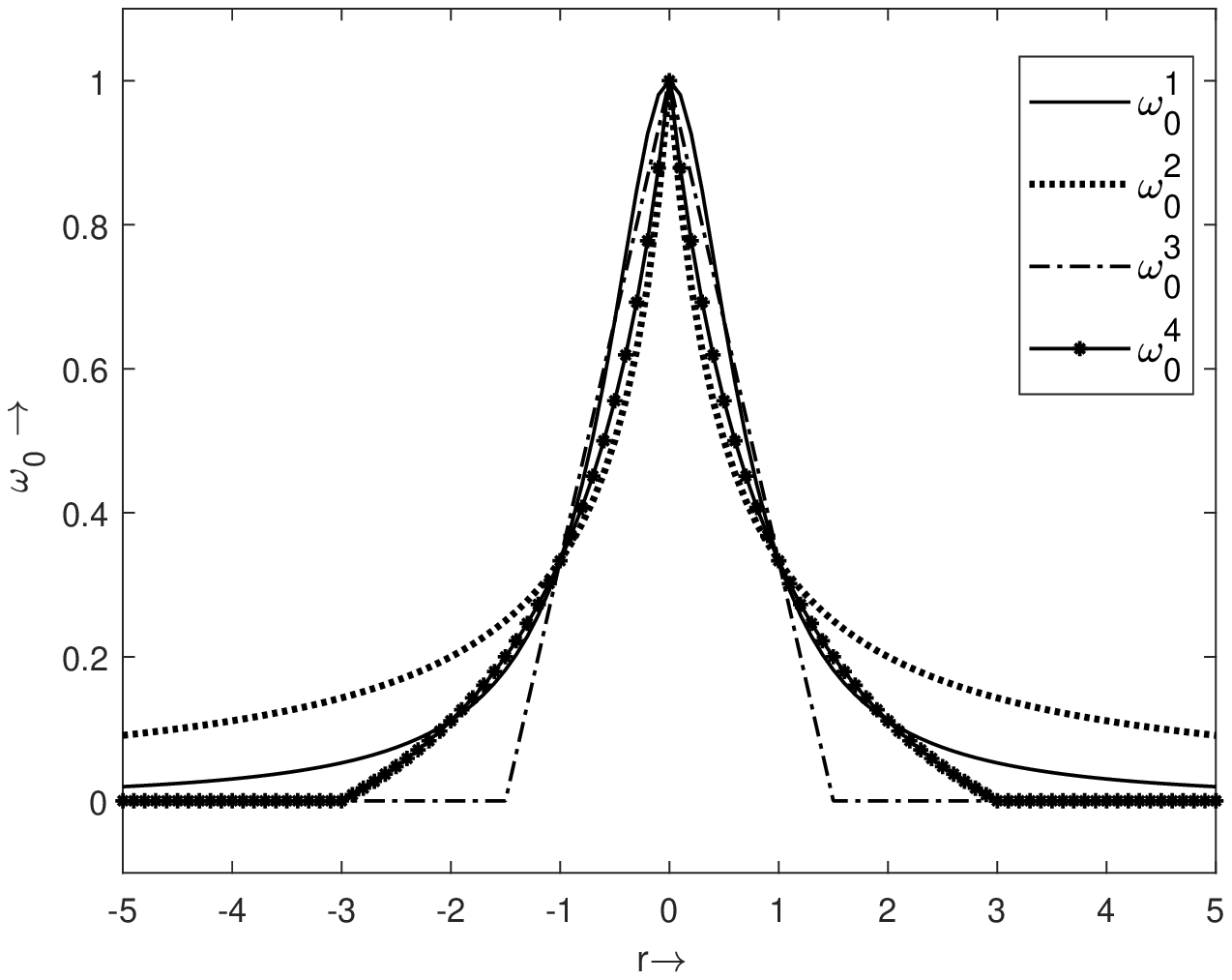}
	\caption{Weight limiter functions $\chi_i$ and corresponding non-linear weights $\omega_0^i= \frac{1}{3}+\frac{2}{3}(1-\chi_i(r)),(i=1,2,3,4)$ }
	\label{fig:limiter1}
\end{figure}

 Based on above observations, the following class of non-linear weights is constructed 
\begin{equation}
\begin{aligned}\label{weights}
\omega_0&=\frac{1}{3}+\frac{2}{3}(1-\chi(r))\\
\omega_1&=1-\omega_0
\end{aligned}
\end{equation}
where $\chi(r)$ is termed as weight limiter function which  must satisfy the following characteristics to yield non-oscillatory third order weights. 
\begin{itemize}
\item[i.] {\bf Non-Oscillatory Conditions:}
\begin{itemize}
	\item[(a)] $\chi(0)=0$ to achieve $\omega_0=1$ and $\omega_1=0$ (Upwind only flux).
	\item[(b)] $\chi(\pm \infty)=\frac{3}{2}$ to get $\omega_0=0$ and $\omega_1=1$ (Centered only flux).	
\end{itemize} 
\item [ii.] {\bf Third order Accuracy Condition:} $\chi$ should be differentiable at $r=1$ and $\chi(1)=1$   to achieve ideal weights $\omega_0=\frac{1}{3}, \; \omega_1=\frac{2}{3}$.
\end{itemize}
The above characterization paves the way to design weight limiter functions such that it achieve third order accuracy even for extrema away from critical points. For example 
\begin{subequations}\label{weightlim}
	\begin{equation}\label{lim1}
	\chi_1(r)=\frac{3r^2}{2r^2+1}
	\end{equation}
	\begin{equation}\label{lim2}
	\chi_2(r) =\frac{3|r|}{2|r|+1}
	\end{equation}
    \begin{equation}\label{lim3}
   	\chi_3(r) = \min(|r|,\frac{3}{2})
    \end{equation}
    \begin{equation}\label{lim4}
    \chi_4(r) = \min(\frac{2|r|}{1+|r|}, \frac{3}{2})
    \end{equation}	
\end{subequations}
 In Figure \ref{fig:limiter1}, weight limiter functions $\chi_i$ \eqref{weightlim} and corresponding weights $\omega_0^i$ \eqref{weights} are given. Another class of limiter with better discontinuity capturing property can be defined as follows
\begin{equation}\label{klimiter}
\chi_5^{k}(r)= \displaystyle \min\left(k|r|,\max\left(1, \frac{3|r|}{2|r|+k}\right)\right)
\end{equation}
where $k\geq 1$ is a constants which guarantee that resulting scheme maintain the third order accuracy in the solution region $1/k\leq |r|\leq k$.
\par
\begin{figure}[h]
 	\centering
 	\includegraphics[scale=0.5]{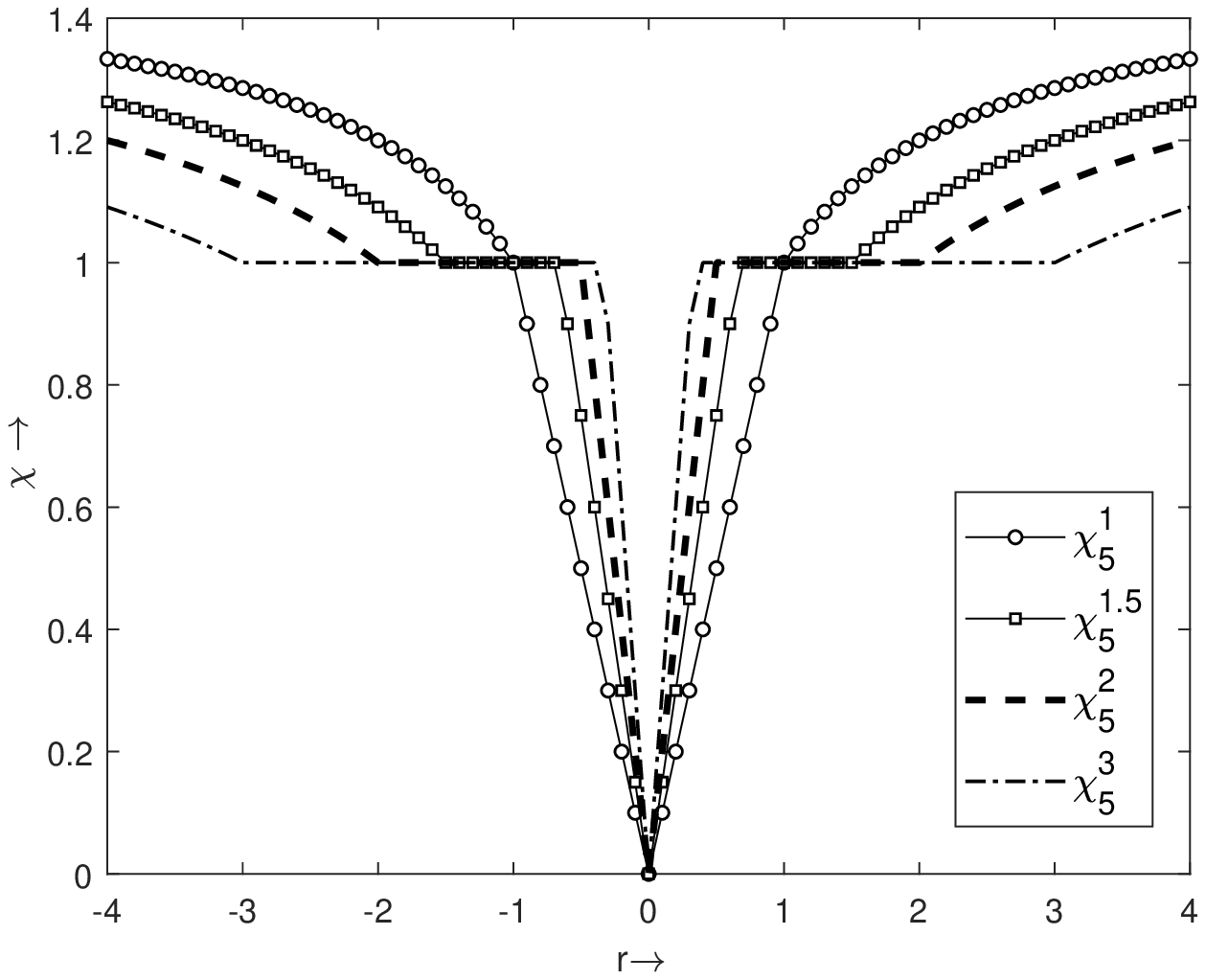}~\includegraphics[scale=0.5]{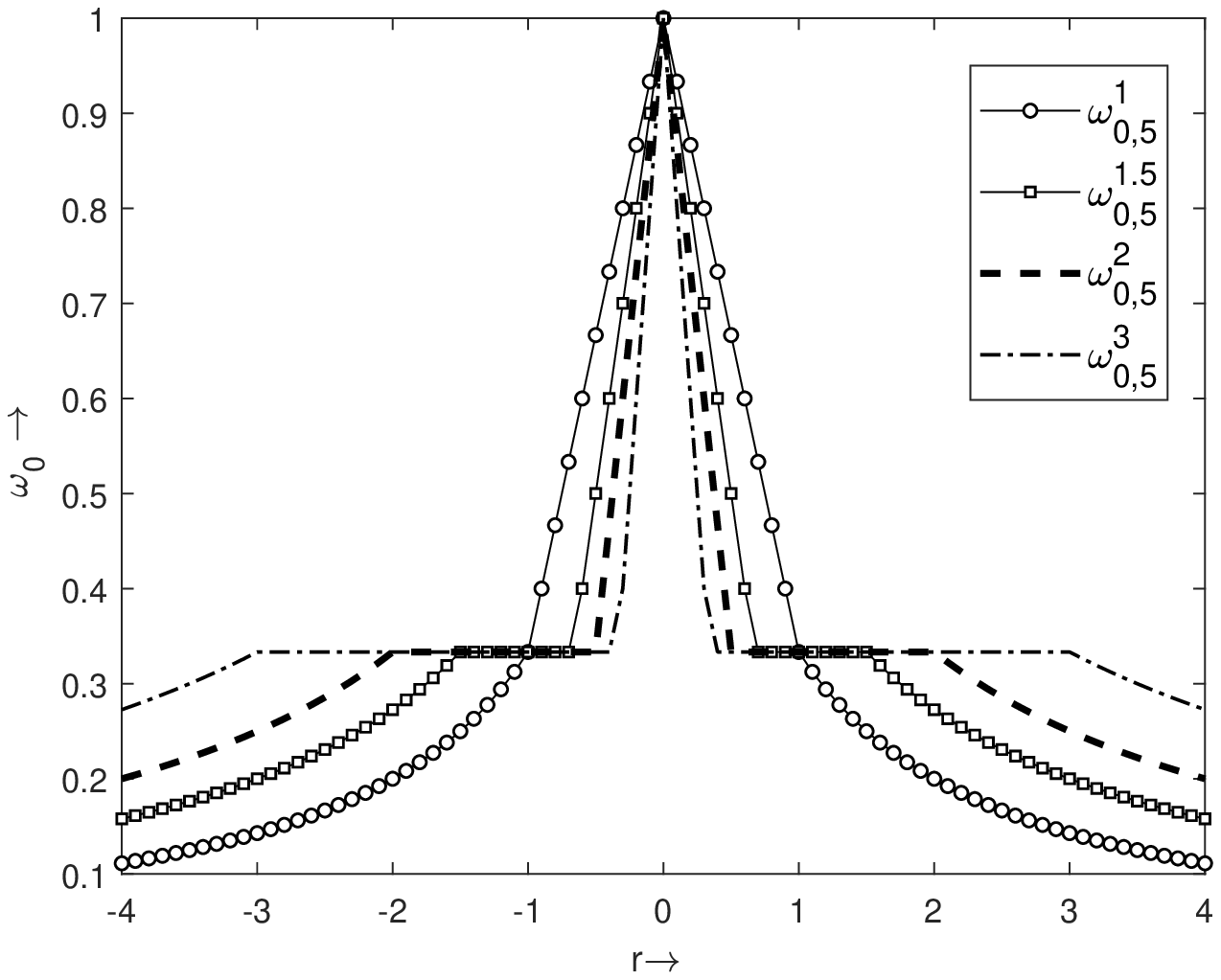}
 	\caption{Weight limiter functions $\chi_5^k$  and corresponding weights $\omega_{0,5}^k= \frac{1}{3}+\frac{2}{3}(1-\chi_5^k(r)),(k=1,1.5,2,3)$ }
 	\label{fig:limiter2}
 \end{figure}

 Note that choice $k>3$ violets the non-oscillatory conditions viz $r\leq 3$ of upwind flux in Lemma \ref{lem3} and $r>\frac{1}{3}$ of centered flux in Lemma \ref{lem2}. This suggests that $k$ must satisfying $1\leq k \leq 3$. Numerical results also support this restriction (see Figure \ref{ptrCondition}(b) and \ref{fig:fig-2c}(c)). Thus it can be concluded weight corresponding to $k=3$ gives maximum possible region for third order accuracy without oscillations. In Figure \ref{fig:limiter2}, weight limiter functions $\chi_{5}^{k}$ \eqref{klimiter} and corresponding weights $\omega_{0,5}^{k}$ \eqref{weights} are given for different values of parameter $k$.
\par Note that all above limiters in \eqref{weightlim} and \eqref{klimiter} are constructed such that they satisfy $\chi(\pm 1)=1$ condition thus guarantee for third order accuracy of scheme at smooth solution region $r\approx 1$ including smooth extrema where $r\approx -1$. 
\begin{remark}
	The weight limiters defined in \eqref{weightlim} and \eqref{klimiter} are functions of smoothness parameter \eqref{smoothnessparameter} similar to the flux limiters used in construction of high resolution shock capturing schemes in \cite{Sweby1984, DUBEY2013325,ZHANG2015114}. Note that, the conditions on flux limiters in \cite{Sweby1984} are derived by imposing the total variation diminishing (TVD) property which demanded that flux limiter functions must vanish for all point of extrema i.e., $r\leq 0$ which cause degenerate accuracy of the TVD schemes at extrema. Numerically this degeneracy leads to clipping error while approximating smooth solution. However, weight limiter function $\chi$ in \eqref{weights} does not necessarily needed to vanish for all $r<0$. Thus weight in \eqref{weights} may retains third order reconstruction even for solution extrema $r\approx -1$ and thus clipping error can be reduce.
\end{remark}
 \begin{remark}
 	The proposed weights using limiter functions \ref{weightlim} or \ref{klimiter} are free of parameters and gives consistent approximations  
 \end{remark}

\section{Numerical results}\label{sec4}
In this section various standard test problems are considered to analyze accuracy and the non-oscillatory behavior of third order WENO3 scheme with proposed weight limiters.  The following name convention is used through out this section.
\begin{itemize}
	\item WENO3-$\omega_0^{i}$ denotes the result obtained by weights \eqref{weights} using function $\chi_i$ from \eqref{weightlim}.
	\item WENO3-$\omega_{0,5}^{l}$ denotes the result obtained by weights \eqref{weights} using function $\chi_5^{k}$ in \eqref{klimiter} for $k=l$.
	\item WENO-$*$ represents the results obtained by various weights (*) in section \ref{sec2} e.g., WENO-$JS3$ represents WENO-JS  weights \eqref{w-js}.
\end{itemize}
\subsection{Linear transport equation}
Consider the linear transport equation,
\begin{equation}\label{eqn4.1}
u_t+a u_x=0,~-1\le x\le 1,~t>0
\end{equation}
with $a=1$ and along with following initial conditions.\\
\subsubsection{Test for Non-oscillatory property} 
\textbf{ Example 1} Consider the  discontinuous initial condition 
\begin{equation}\label{Ltest2a}
u_0(x)=\left\{\begin{array}{cc}
1 & |x|\leq 0.3,\\
0 & else.
\end{array}\right.
\end{equation}
\textbf{ Example 2} Consider the smooth initial condition with sharp turn 
\begin{equation}\label{Ltest2b}
u_0(x)=\left\{\begin{array}{cc} [0.5+0.5\cos(\omega(x-x_c))]^4 & if |x-x_c|<\sigma \\
0 & otherwise \end{array}\right\}
\end{equation}
where computational domain is $0\leq x\leq 1, \omega=5\pi$, $x_c=0.5$ and $\sigma=0.2$.\\

\begin{figure}[htb!] 
	\begin{tabular}{cc}
		\hspace{-1.2cm}
		\includegraphics[scale=0.55]{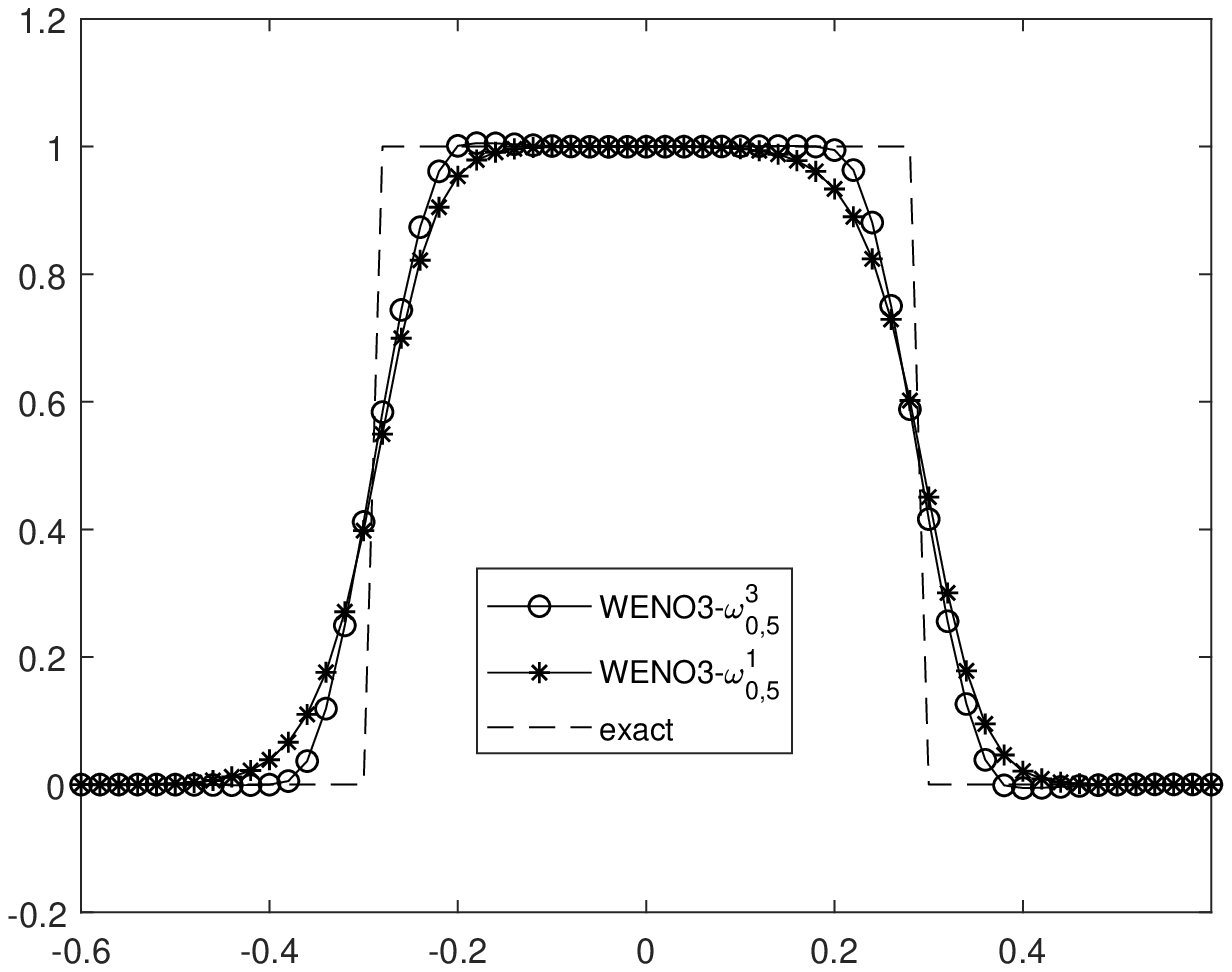} &
		\includegraphics[scale=0.55]{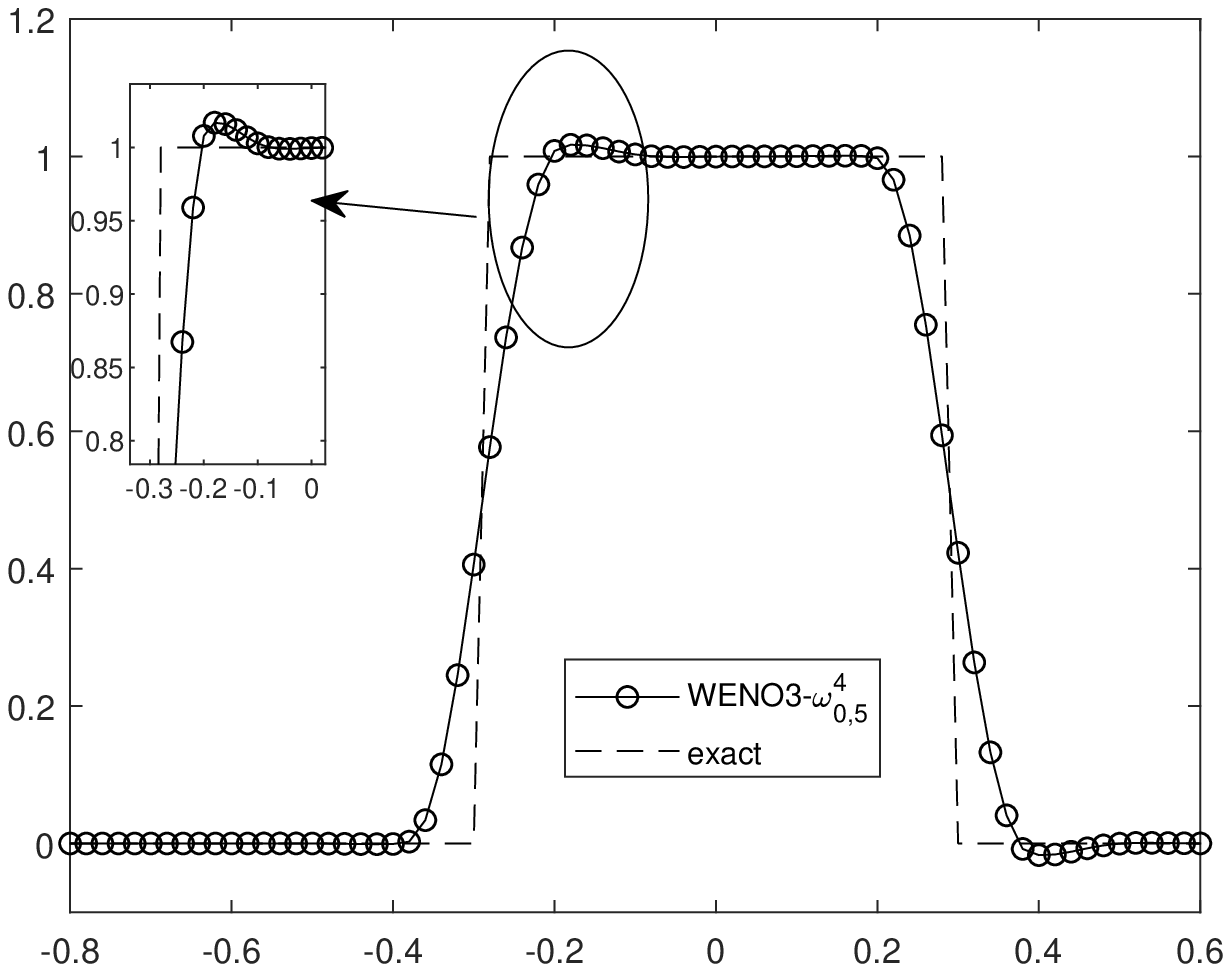} \\
		\textbf{a}&\textbf{b}
	\end{tabular}
	\caption{Solution of linear equation \eqref{eqn4.1} of WENO3-$\omega_{0,5}^{k}$ with initial condition \eqref{Ltest2a} at $t=2$ with square initial condition $\frac{\Delta t}{\Delta x} =0.5$: Effect of $k$ on oscillatory behavior.}	\label{ptrCondition}	
\end{figure}

In Figure \ref{ptrCondition}, numerical results for example 1 are given using weight limiter \eqref{klimiter} for $k=1, 3$. It can be seen in Figure \ref{ptrCondition}\textbf{(a)} that numerical approximation with both choices of $k$ is non-oscillatory however choice $k=1$ yields more diffusive solution compared to the choice $k=3$. Figure \ref{ptrCondition}\textbf{(b)} shows oscillatory behavior of weight limiter \eqref{klimiter} for the choice $k=4$. We mention that the numerical solution obtained by weight limiter functions \eqref{lim1}-\eqref{lim4} are not shown as they falls between these two approximations.

\begin{figure}[htb!] 
	\begin{tabular}{cc}
		\hspace{-1.2cm}
		\includegraphics[scale=0.55]{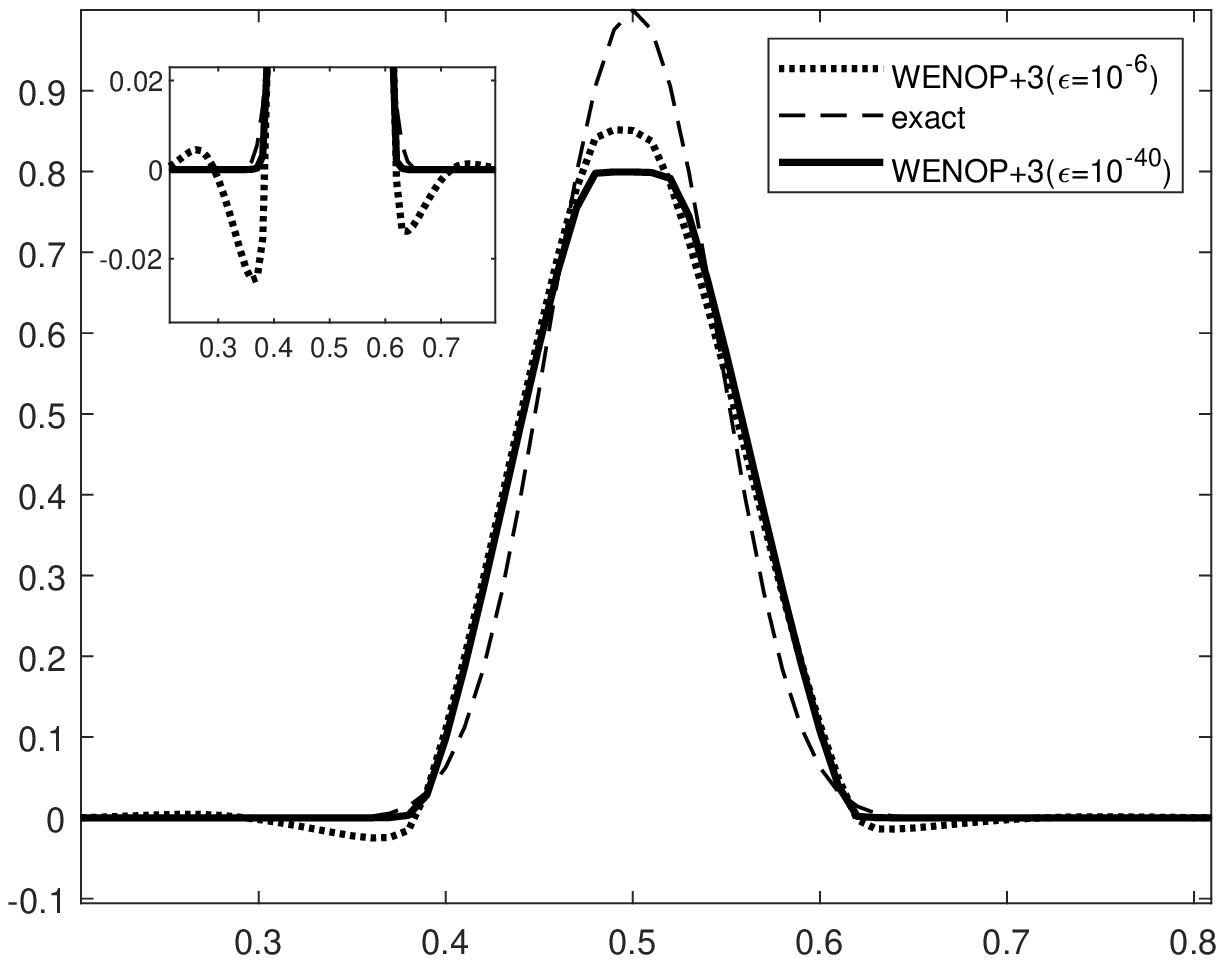} &
		\includegraphics[scale=0.55]{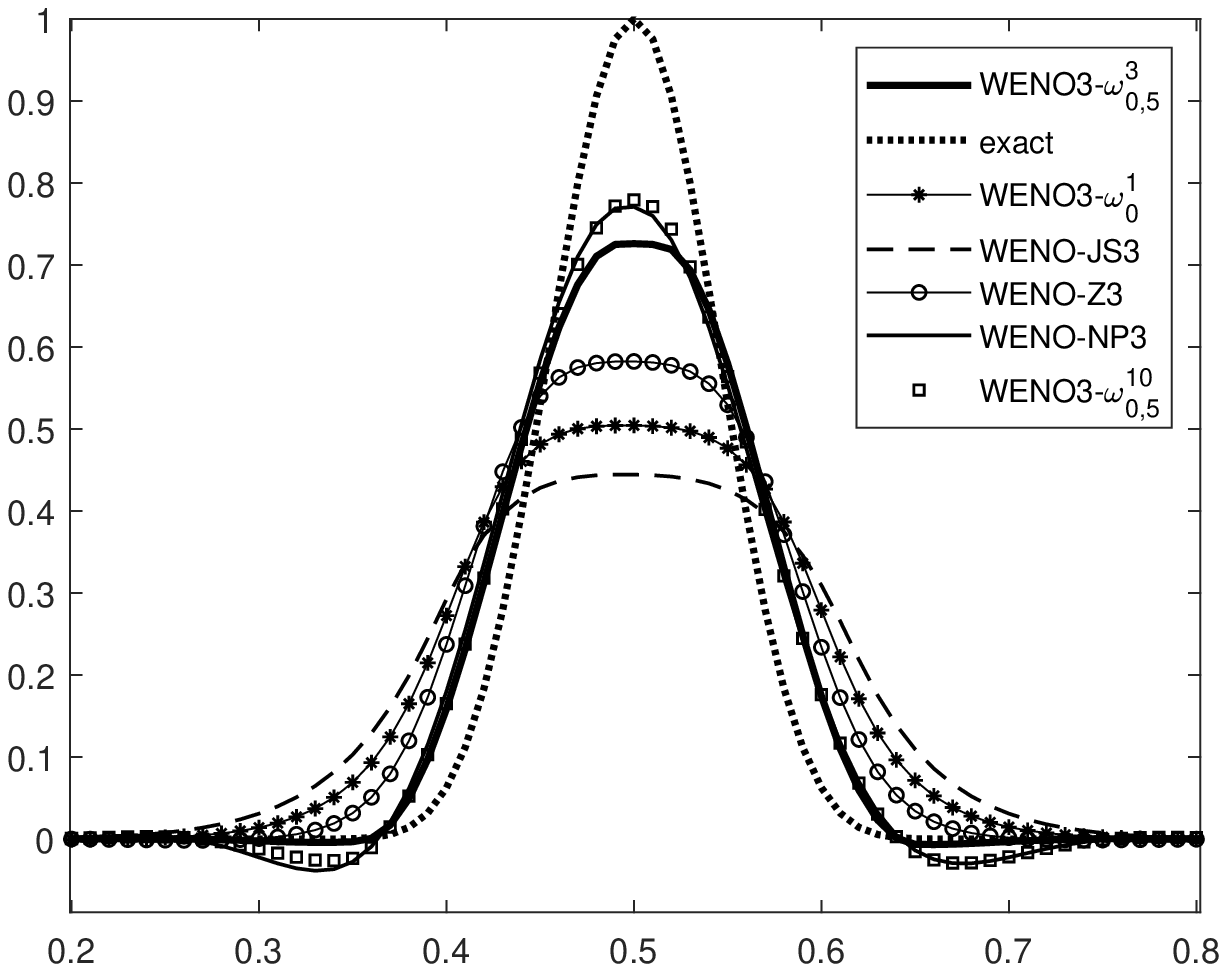} \\
		\textbf{a}&\textbf{b}
	\end{tabular}
	\caption{Solution for initial condition \eqref{Ltest2b} at $t=10$ with square initial condition $\frac{\Delta t}{\Delta x} =0.5$ (a) Oscillatory solution by WENO-P+3 for parameter $\epsilon =1e-6$ (b) Comparison of solution}	
	\label{FigCompare}	
\end{figure}

In Figure \ref{FigCompare}, numerical solution by WENO3 scheme using various existing and newly proposed weights are given and compared. It can be noted that proposed WENO-$\omega_{0,5}^k$ weights gives better approximation to maxima without oscillations for $k=3$. The results with choice $k=10$ gives improved approximation for smooth maxima but exhibits small undershoot in bottom smooth region between $0.3\leq x\leq 0.4$ and $0.6\leq x\leq 0.7$. This behavior is also true for WENO-$NP3$ and WENO-P+3 for the choice $\epsilon =1e-6$ in \eqref{WENOp3}.  The result obtained by WENO-P+3 for the choice $\epsilon =1e-40$ in  \eqref{WENOp3} does not show any such oscillations though there is a drop in the peak of smooth maxima with flatness compared to the choice $\epsilon =1e-6$. These results show that in order to achieve improved approximation without clipping error for smooth maxima one has to compromise on the non-oscillatory property.
\subsubsection{Accuracy test}\label{sec4b}
Consider equation \eqref{eqn4.1} with smooth initial conditions\\
\textbf{ Example 3} 
\begin{equation}\label{eqn4.2}
u_0(x)=sin(\pi x)~\text{in the periodic domain [-1,1]}
\end{equation}
\textbf{ Example 4} 
\begin{equation}\label{eqn4.3}
u_0(x)=sin^4(\pi x)~\text{in the periodic domain [0,1]}
\end{equation}
\begin{figure}[htb!] 
	\begin{tabular}{cc}
		\centering 
		\includegraphics[scale=0.5]{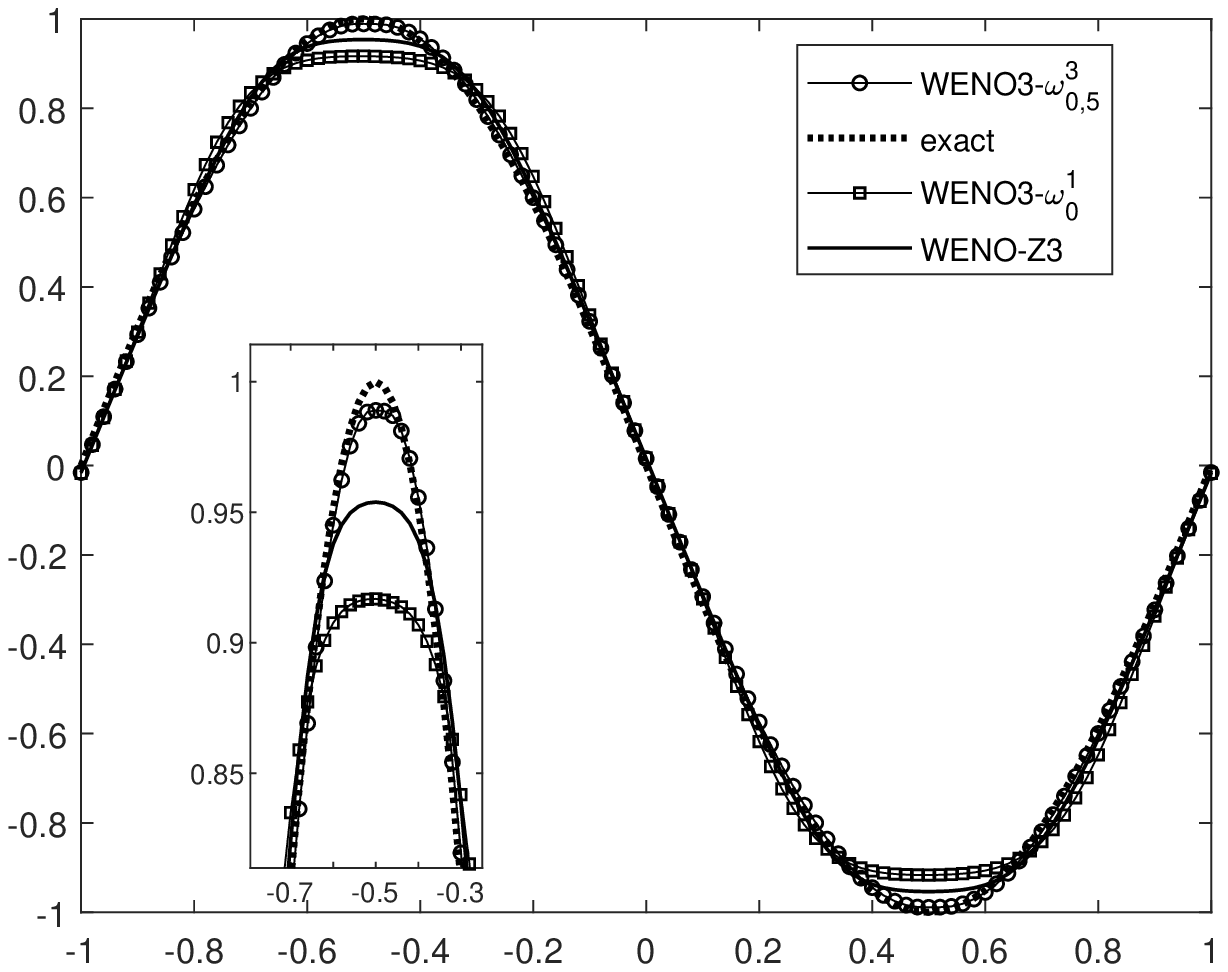} &
		\includegraphics[scale=0.5]{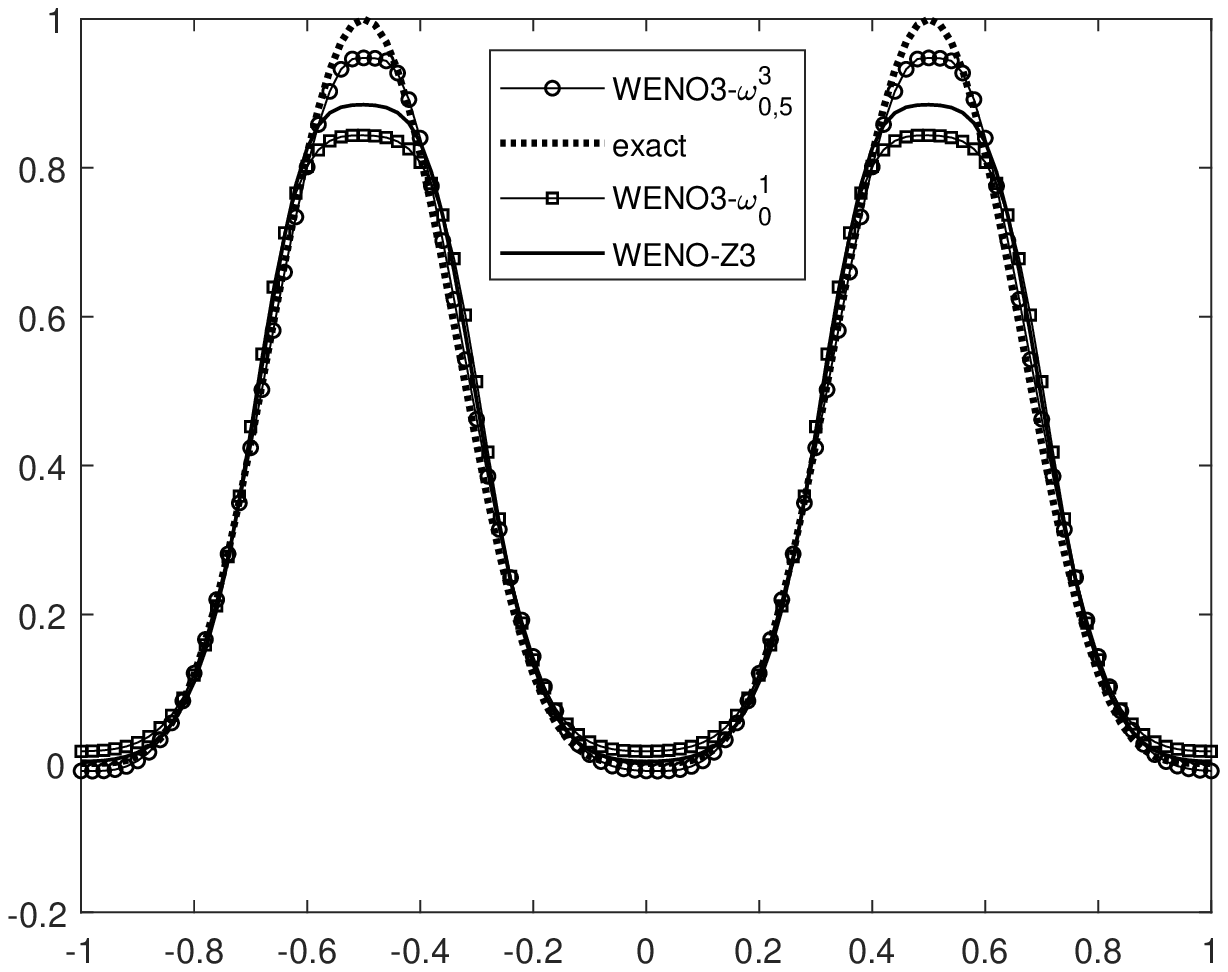} \\
		\textbf{a}&\textbf{b}
	\end{tabular}
	\caption{Comparison of Solution of linear equation \eqref{eqn4.1} of WENO-Z3, WENO3-$\omega_0^1$, WENO3-$\omega_{0,5}^{3}$ with smooth \textbf{(a)}initial condition \eqref{eqn4.2} $t=30$, $\frac{\Delta t}{\Delta x} =0.5$, \textbf{(b)} initial condition \eqref{eqn4.3} $t=4$ $\frac{\Delta t}{\Delta x} =0.25$.}
	\label{clipping}	
\end{figure}
In Figure \ref{clipping}, numerical approximation using proposed weights \eqref{lim1} and \eqref{klimiter} for smooth initial conditions is given and compared with weight WENO-$Z3$ \cite{DON2013347} .
It can be seen from results that proposed weight $\omega_{0,5}^3$ gives significantly better approximation for smooth extrema compared to the WENO-$Z$.\\ 

In Tables \ref{tab:1} and \ref{tab:2}, the convergence rate in both $L^1$ and $L^{\infty}$ error is given for the WENO3 scheme using proposed weights $\omega_0^1$ and $\omega_{0,5}^1$ for the problem \eqref{eqn4.1} with initial condition \eqref{eqn4.2}. Similarly Tables  \ref{tab:3} and \ref{tab:4} shows the convergence rate of these weights for the problem \eqref{eqn4.1} with initial condition \eqref{eqn4.3}. We remark that, computationally weights corresponding to all the limiter functions in \eqref{weightlim} and \eqref{klimiter} give more or less  similar third order accuracy however weight limiters \eqref{lim1}, \eqref{lim2} and \eqref{klimiter} for $k=1$ give better consistent third order convergence rate due to their smooth nature. 
\begin{table}[htb!]
	
	\centering
	\begin{tabular}{|c|c|c|c|c|}
		\hline N  & WENO3-$\omega_0^1$  &  Rate &    WENO3-$\omega_0^1$    & Rate \\ 
		& $L^\infty$ error  &   &    $L^1$ error    &  \\ 
		\hline 80 & 1.24987e-02 &  -Inf & 4.44919e-03 & -Inf \\ 
		\hline 160 & 3.81984e-03 &  1.71 & 8.45946e-04 & 2.39 \\ 
		\hline 320 & 5.91893e-04 &  2.69 & 9.35838e-05 & 3.18 \\ 
		\hline 640 & 7.22438e-05 &  3.03 & 8.53235e-06 & 3.46 \\ 
		\hline 1280 & 8.88347e-06 &  3.02 & 6.51676e-07 & 3.71 \\ 
		\hline 2560 & 9.28233e-07 &  3.26 & 4.18069e-08 & 3.96 \\ 
		\hline   
	\end{tabular} 
	\caption{Rate of convergence for IC \eqref{eqn4.2} at final time $t=0.5$ with $CFL=0.25.$}
	\label{tab:1}
\end{table}
\begin{table}[htb!]
	
	\centering
	\begin{tabular}{|c|c|c|c|c|}
		\hline N  & WENO3-$\omega_{0,5}^1$  &  Rate &    WENO3-$\omega_{0,5}^1$    & Rate \\ 
		& $L^\infty$ error  &   &    $L^1$ error    &  \\  
		\hline 80 & 8.83849e-03 &  -Inf & 2.68444e-03 & -Inf \\ 
		\hline 160 & 2.53953e-03 &  1.80 & 4.69982e-04 & 2.51 \\ 
		\hline 320 & 3.37745e-04 &  2.91 & 4.75547e-05 & 3.30 \\ 
		\hline 640 & 3.65303e-05 &  3.21 & 4.15503e-06 & 3.52 \\ 
		\hline 1280 & 4.45903e-06 &  3.03 & 3.16314e-07 & 3.72 \\ 
		\hline 2560 & 4.65153e-07 &  3.26 & 2.03236e-08 & 3.96 \\ 
		\hline 
	\end{tabular} 
	\caption{Rate of convergence for IC \eqref{eqn4.2} at final time $t=0.5$ with $CFL=0.25$.}
	\label{tab:2}
\end{table}	 
\begin{table}[htb!]
	\centering
	\begin{tabular}{|c|c|c|c|c|}
		\hline N  & WENO3-$\omega_0^1$  &  Rate &    WENO3-$\omega_0^1$    & Rate \\ 
		& $L^\infty$ error  &   &    $L^1$ error    &  \\  
	\hline 80 & 4.76670e-02 &  -Inf & 1.85105e-02 & -Inf \\ 
	\hline 160 & 1.83789e-02 &  1.37 & 4.16305e-03 & 2.15 \\ 
	\hline 320 & 5.61142e-03 &  1.71 & 7.71866e-04 & 2.43 \\ 
	\hline 640 & 9.96515e-04 &  2.49 & 9.43049e-05 & 3.03 \\ 
	\hline 1280 & 1.14091e-04 &  3.13 & 8.67465e-06 & 3.44 \\ 
	\hline 2560 & 1.41383e-05 &  3.01 & 7.30714e-07 & 3.57 \\ 
	\hline
	\end{tabular} 
	\caption{Rate of convergence for IC \eqref{eqn4.3} at final time $t=0.5$ with $CFL=0.25$.}
	\label{tab:3}
\end{table}
\begin{table}[htb!]
	
	\centering
	\begin{tabular}{|c|c|c|c|c|}
		\hline N  & WENO3-$\omega_{0,5}^1$  &  Rate &    WENO3-$\omega_{0,5}^1$    & Rate \\ 
		& $L^\infty$ error  &   &    $L^1$ error    &  \\
		\hline 80 & 3.46967e-02 &  -Inf & 1.07368e-02 & -Inf \\ 
		\hline 160 & 1.29620e-02 &  1.42 & 2.39369e-03 & 2.17 \\ 
		\hline 320 & 3.65421e-03 &  1.83 & 4.12535e-04 & 2.54 \\ 
		\hline 640 & 5.80597e-04 &  2.65 & 4.72786e-05 & 3.13 \\ 
		\hline 1280 & 5.80477e-05 &  3.32 & 4.20458e-06 & 3.49 \\ 
		\hline 2560 & 7.12878e-06 &  3.03 & 3.54383e-07 & 3.57 \\ 
		\hline 
	\end{tabular} 
	\caption{Rate of convergence for IC \eqref{eqn4.3} at final time $t=0.5$ with $CFL=0.25.$}
	\label{tab:4}
\end{table}

\subsection{Inviscid Burgers' equation}
\begin{figure}[htb!]
	\begin{tabular}{ccc}
	\hspace{-1cm}\centering
	\includegraphics[scale=0.35]{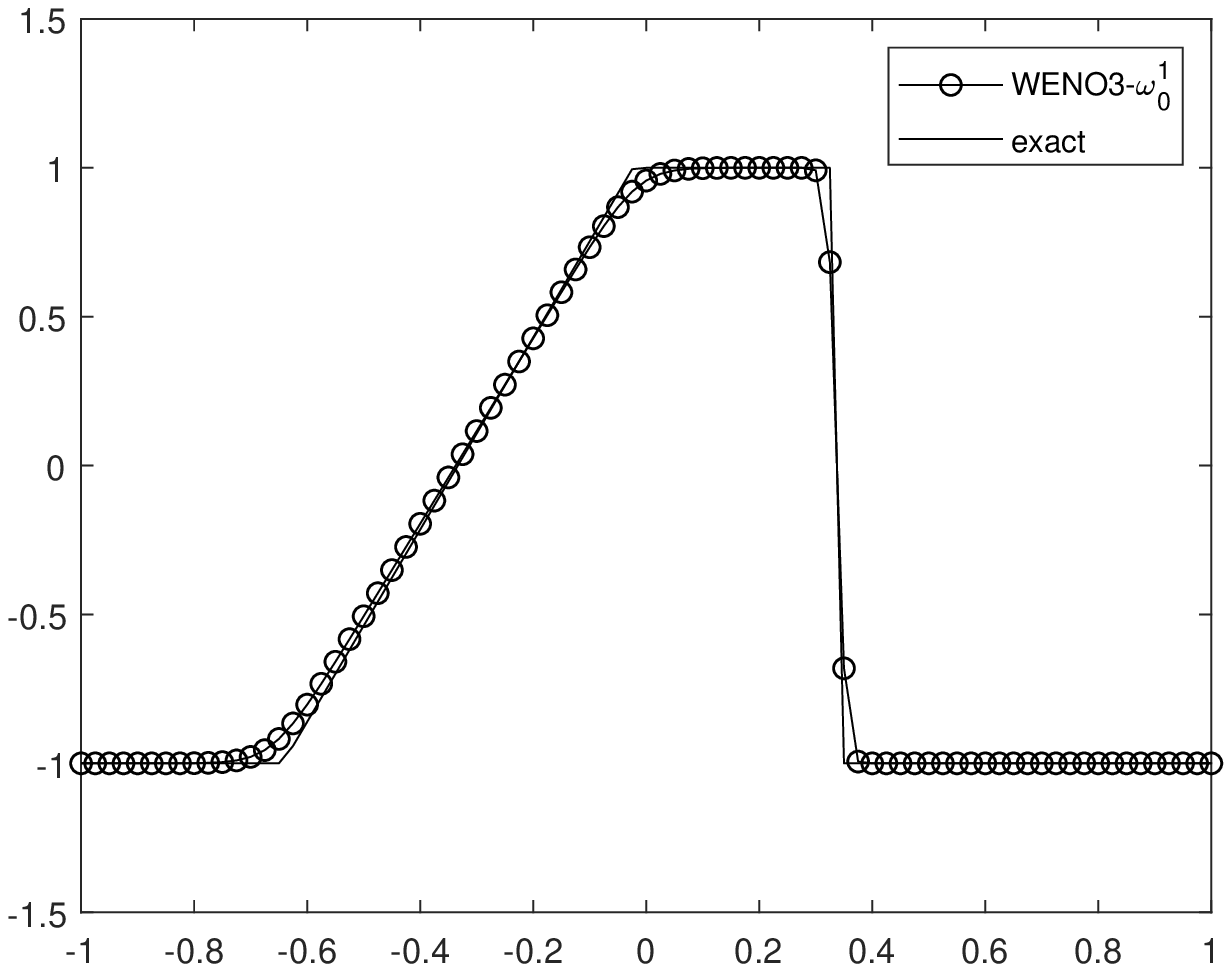} &
	\includegraphics[scale=0.35]{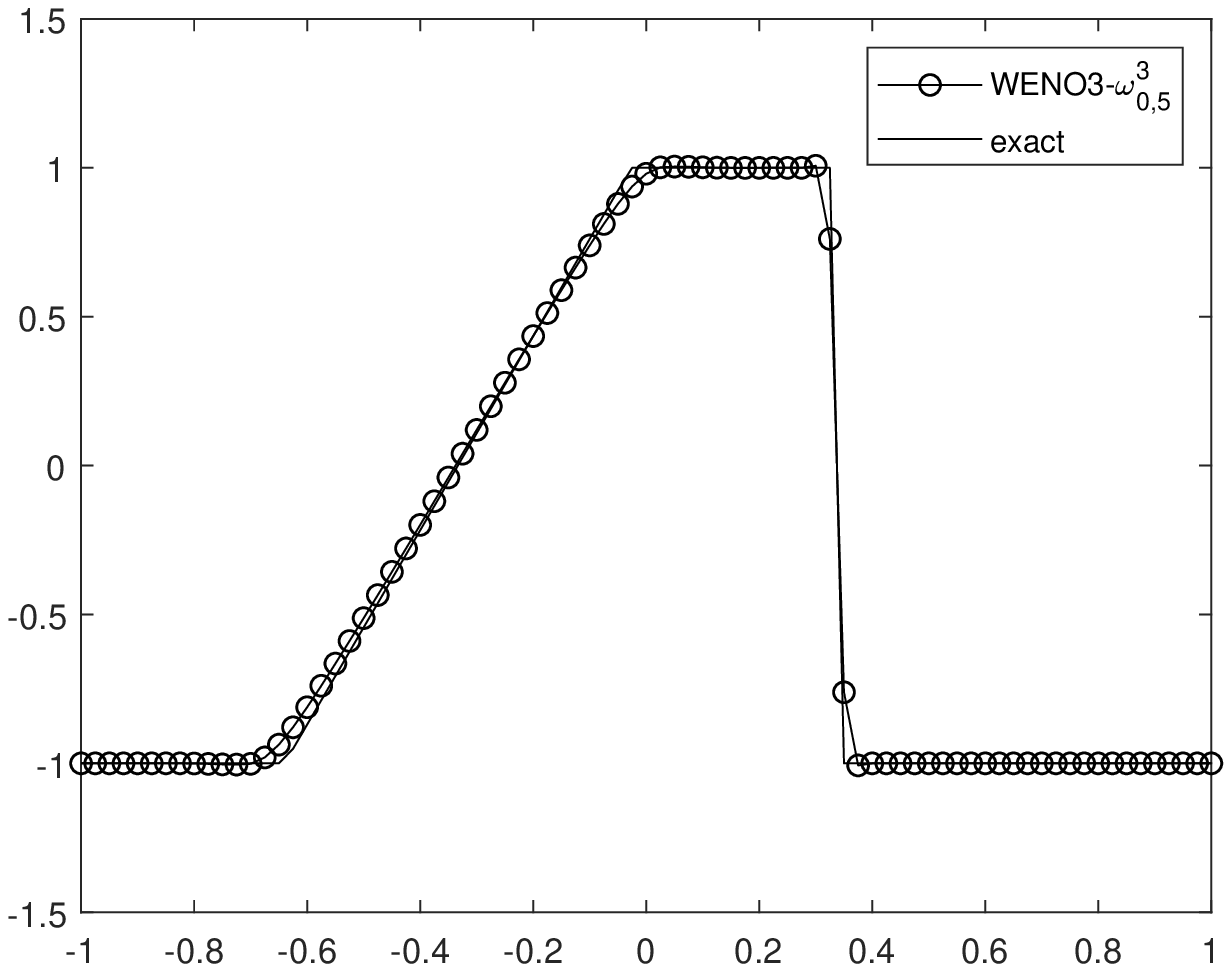}& \includegraphics[scale=0.35]{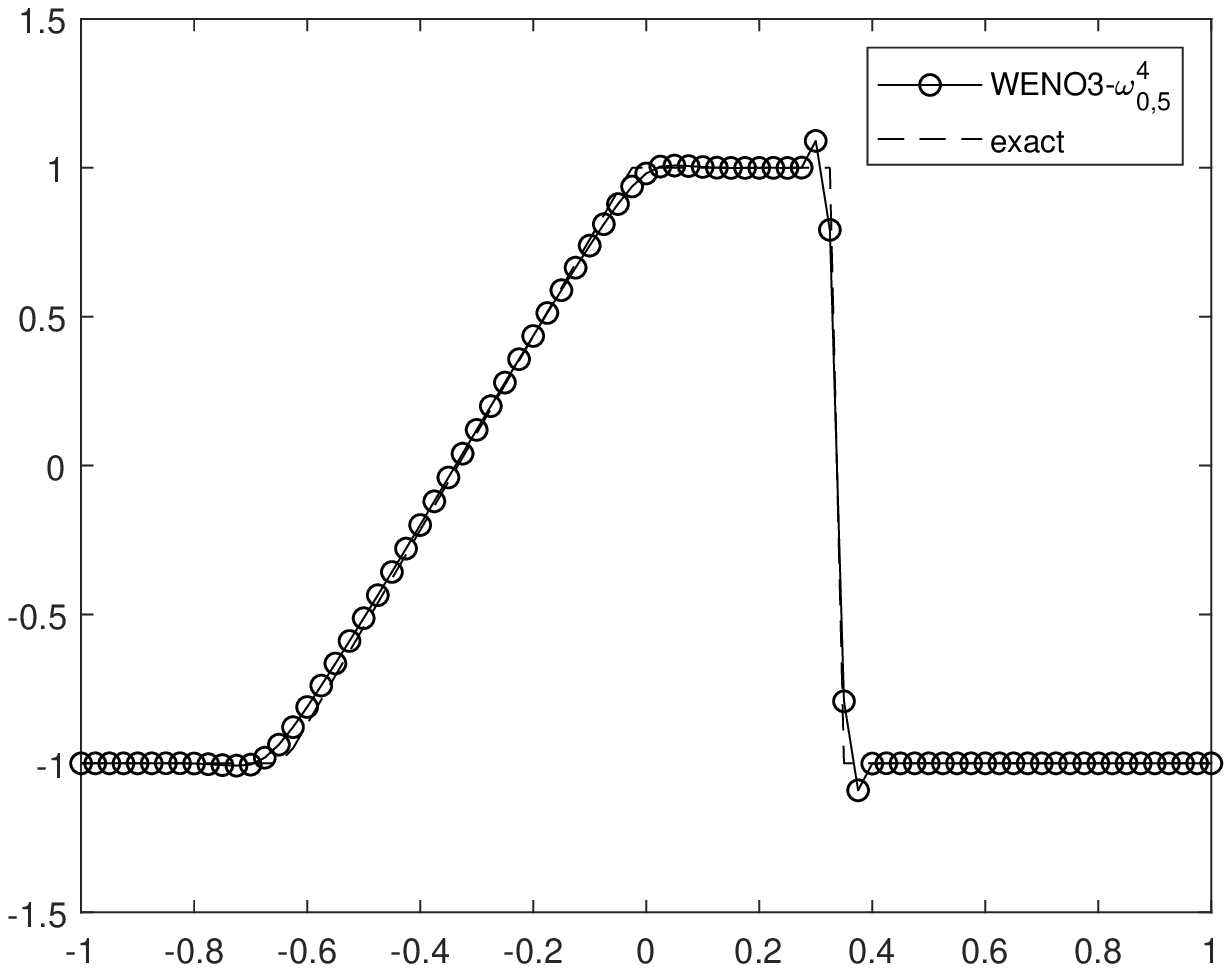}\\
	\textbf{a}&\textbf{b} & \textbf{c}
	\end{tabular}
	\caption{Solution of Burgers' equation with initial condition \eqref{eqn5.2} at $t=0.3$ with $CFL=0.5$ in the periodic domain of $[-1,1]$ with $100$ grids.}
	\label{fig:fig-2c}
\end{figure}

Consider the Burger equation
\begin{equation}\label{eqn5.1}
u_t+\left(\frac{u^2}{2}\right)_x=0,~-1\le x\le 1,~t>0
\end{equation}
with discontinuous initial condition
\begin{equation}\label{eqn5.2}
u(x,0)=
\begin{cases}
1~~ \text{if}~~ |x|<\frac{1}{3}\\
-1~ \text{if}~~ \frac{1}{3}\le |x|\le 1.
\end{cases}
\end{equation}
The exact solution of Burger equation corresponding to initial condition \eqref{eqn5.2} consists a left rarefaction and a steady shock at point $x=\frac{1}{3}$. Figure \ref{fig:fig-2c}-\textbf{(a)} and \textbf{(b)} shows WENO3 scheme using weight limiter $\chi_5^{k},(k=3)$ capture the head and tail of the rarefaction more accurately compared to using $\chi_1$. Also for both choices of weights, resolution of steady shock is similar and with no oscillations. Figure \ref{fig:fig-2c}-\textbf{(c)} shows the oscillatory behavior of weight limiter \eqref{klimiter} $\chi_5^{k}$ for choice $k=4$.\\
\subsection{1D Euler equations}
The numerical simulations are performed on the one-dimensional Euler equations which are given by
\begin{equation}\label{euler}
\left(\begin{array}{c}
\rho\\
\rho u\\
E  \end{array}\right)_t+\left(\begin{array}{c}
\rho u\\
\rho u^2+p\\
u(E+p)  \end{array}\right)_x=0
\end{equation}
where $\rho,~u,~E,~p$ are the density, velocity, total energy and pressure respectively. The system \eqref{euler} represents the conservation of mass, momentum and energy. The total energy for an ideal polytropic gas is defined as 
\begin{equation}\nonumber
E=\frac{p}{\gamma -1}+\frac{1}{2}\rho u^2,
\end{equation}
and the eigen values of the Jacobian matrix $A(U)=\frac{\partial F}{\partial U}$ are 
\begin{equation}\label{eqn5.5}
\lambda_1(u)=u-c,~~\lambda_2(u)=u,~~\lambda_3(u)=u+c,
\end{equation}
where $U=\left(\begin{array}{c}
\rho\\
\rho u\\
E  \end{array}\right),~~F=\left(\begin{array}{c}
\rho u\\
\rho u^2+p\\
u(E+p)  \end{array}\right)$ and $\gamma$ is the ratio of specific heats and its values is taken as \eqref{eqn5.5}.
\subsubsection{Sod's shock tube problem }
\begin{figure}[htb!]
	\begin{tabular}{cc}
		\hspace{-1.5cm}
		\includegraphics[scale=0.55]{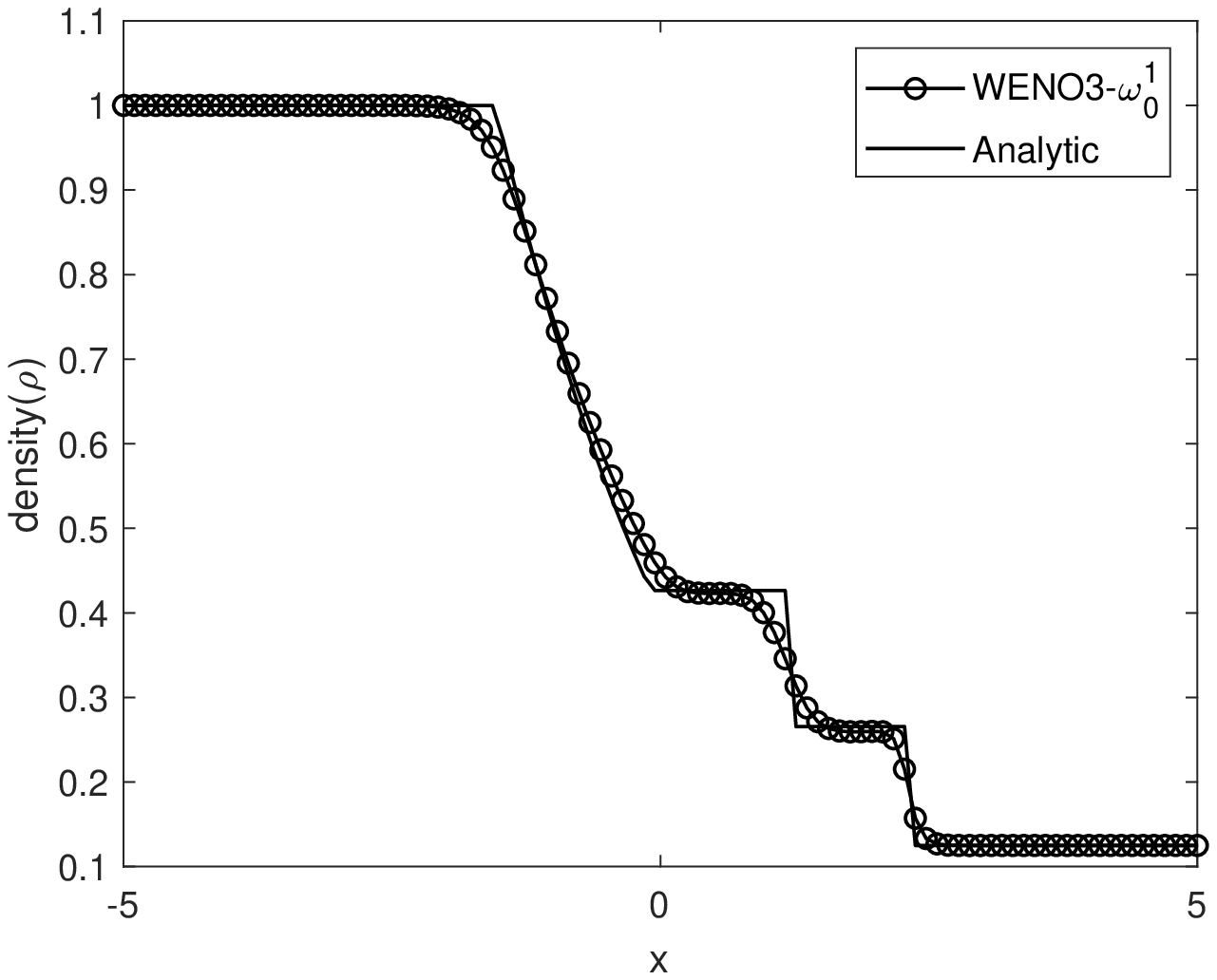} &
		\includegraphics[scale=0.55]{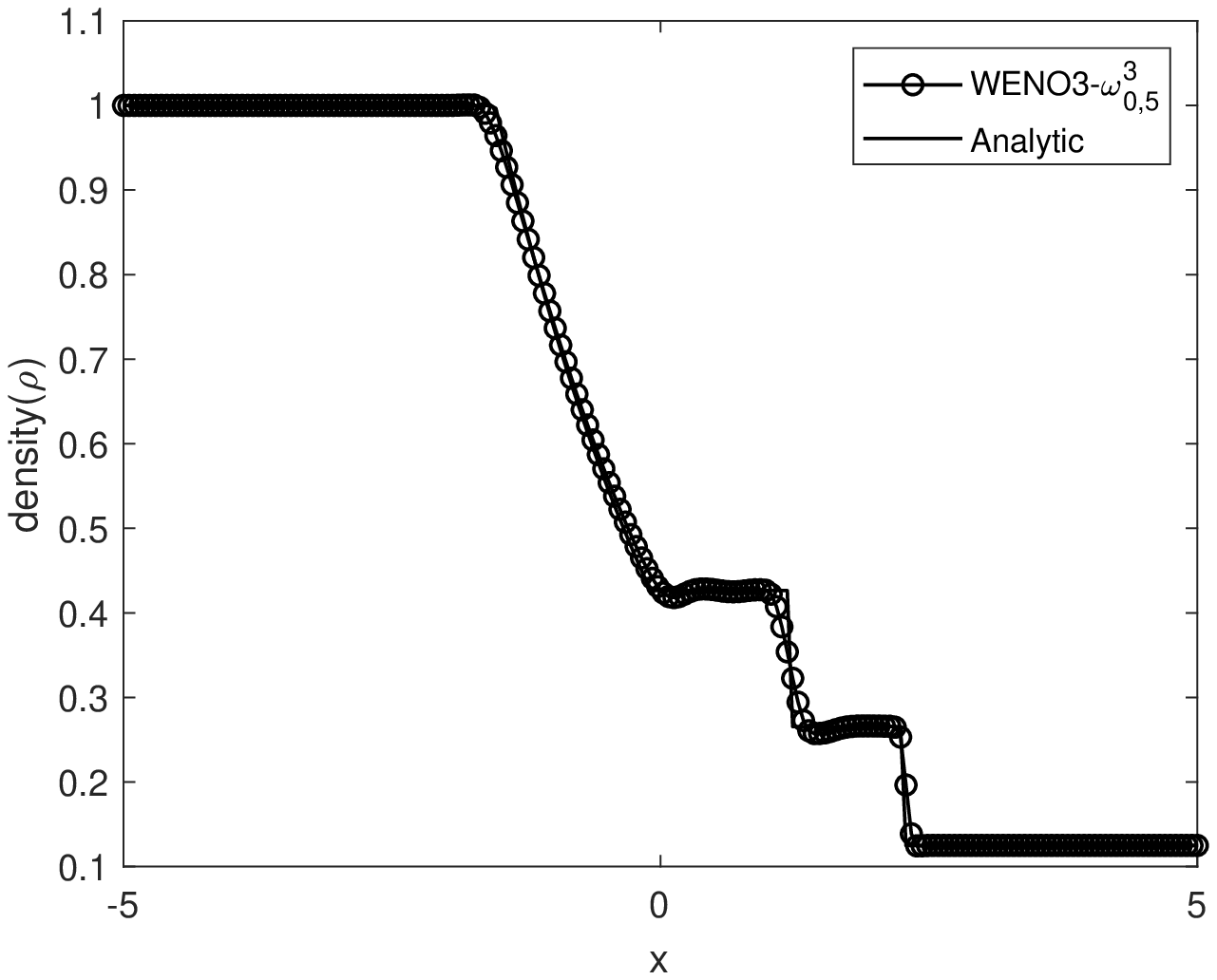}\\
		\textbf{a}&\textbf{b}
	\end{tabular}
	\caption{Solution of Sod's shock tube test\eqref{sod} using weight \textbf{(a)} WENO3-$\omega_0^1$ \textbf{(b)}  WENO3-$\omega_{0,5}^3$. The computational domain $[-5,5]$ is partitioned with $200$ grids and results are computed using $CFL=0.4$ at $t=1.3$.}
	\label{fig:fig-3d}
\end{figure}

Sod's shock tube problem \cite{sod1978survey} is one-dimensional Euler system \eqref{euler} with Riemann data
\begin{equation}\label{sod}
(\rho,~u,~p)(x,0)=
\begin{cases}
(1,~0,~1)~~~~~~~~-5\le x<0,\\
(0.125,~0,~0.1)~~~~0\le x\le 5,
\end{cases}
\end{equation}
in the computational domain $[-5,5]$ at final time $t=1.3$. The density profile of the solution consists of a rarefaction region as well as a shock and contact discontinuity region. In Figure \ref{fig:fig-3d}(a) and \ref{fig:fig-3d}(b) results are given using weights $\omega_0^1$ and $\omega_{0,5}^3$ respectively. It can be seen that WENO3-$\omega_0^1$ captures the shock and contact discontinuities correctly without oscillations whereas WENO3-$\omega_{0,5}^3$ shows a small local undershoot at the foot of the rarefaction. It is needed to be mentioned here that  this undershoot does not grow on further refinement of the mesh and is similar to the results obtained by WENO3-$YC$ scheme.    
\subsubsection{Lax's shock tube problem }
\begin{figure}[htb!]
	\begin{tabular}{cc}
		\hspace{-1.5cm}
		\includegraphics[scale=0.55]{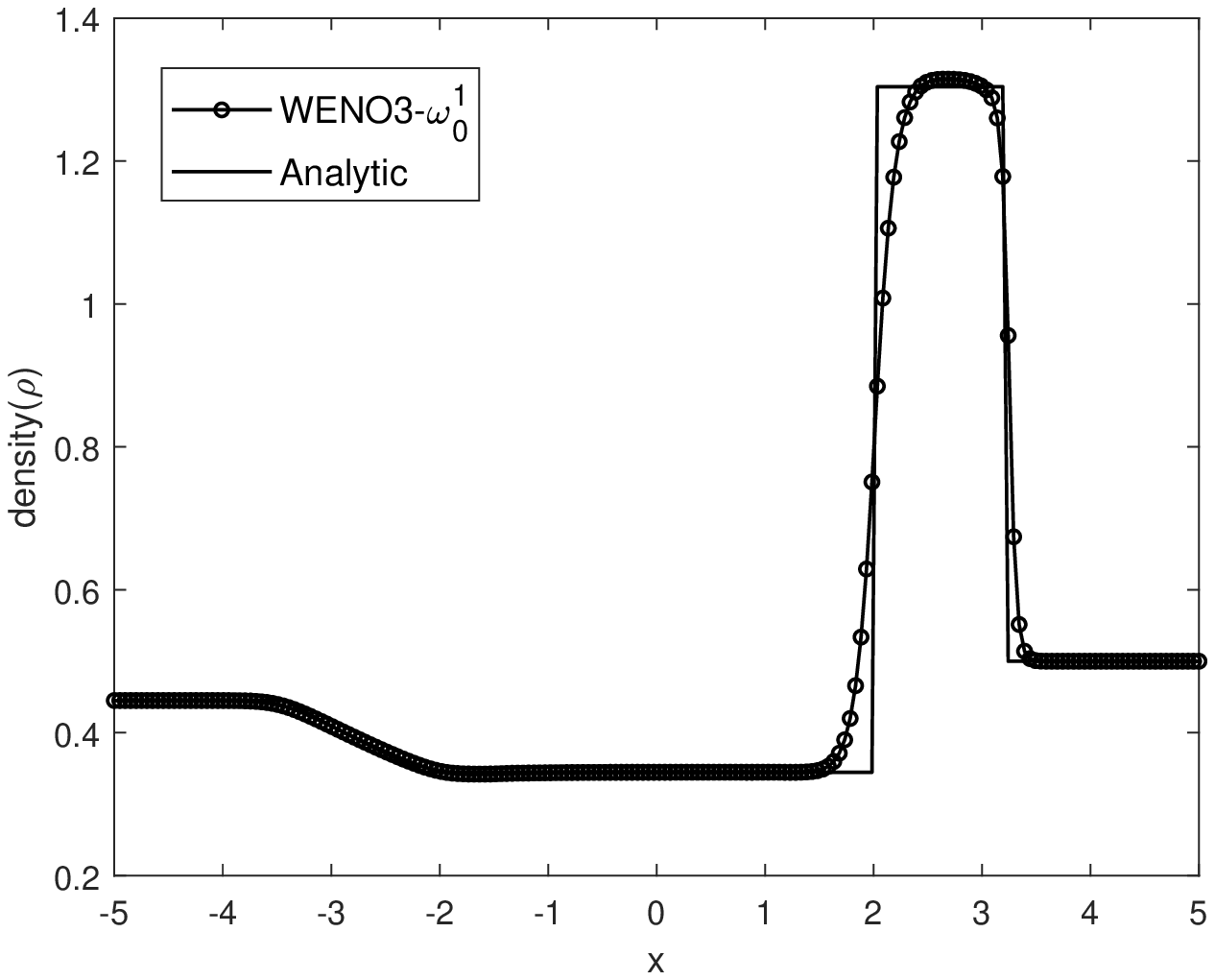} &
		\includegraphics[scale=0.55]{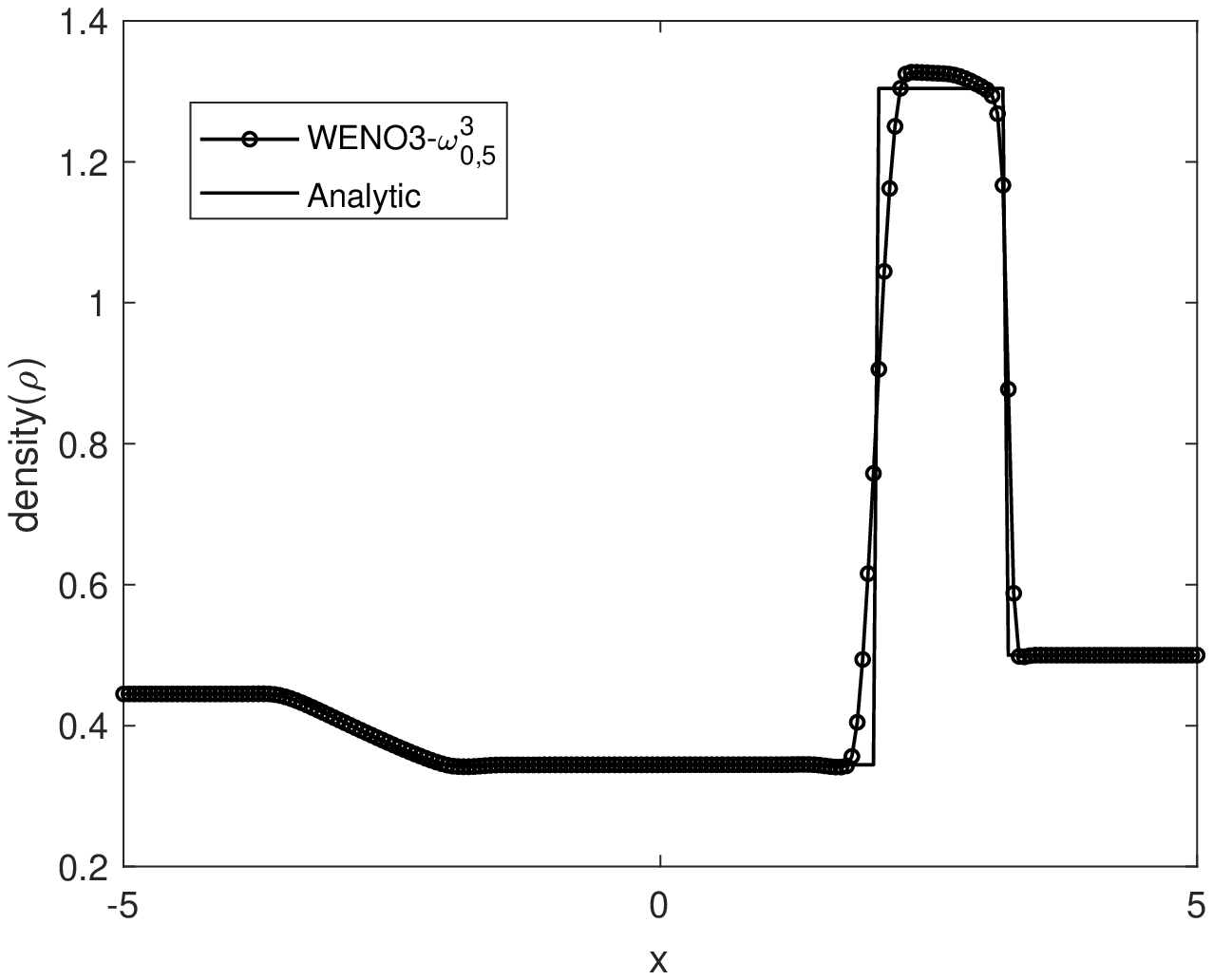}\\
		\textbf{a}&\textbf{b}
	\end{tabular}
	\caption{Solution of Lax's shock tube problem \eqref{lax} using weight \textbf{(a)} WENO3-$\omega_{0}^1$ and \textbf{(b)} WENO3-$\omega_{0,5}^3$. Solution is given in the domain $[-5,5]$ with $200$ grids with $CFL=0.25$ at $t=1.3$.}
	\label{fig:fig-4c}
\end{figure}

Lax's shock tube problem \cite{lax1954weak} with initial condition
\begin{equation}\label{lax}
(\rho,~u,~p)(x,0)=
\begin{cases}
(0.445,~0.698,~3.528)~~~~~-5\le x<0,\\
(0.5,~0,~0.571)~~~~~~~~~~~~~~~~~0\le x\le 5,
\end{cases}
\end{equation}
is an one-dimensional Euler system of equation \eqref{euler} in the computational domain $[-5,5]$ and is run upto $t=1.3$ with zero gradient boundary conditions. The numerical results of density profiles along with reference solutions are displayed in Figure \ref{fig:fig-4c}(a) and \ref{fig:fig-4c}(b). In this test WENO3-$\omega_{0,5}^3$ resolves the discontinuities more crisply compared to WENO3-$\omega_{0}^1$ with spurious oscillations.
\subsubsection{Shu-Osher test }
\begin{figure}[htb!]
	\begin{tabular}{cc}
		\hspace{-1.0cm}
		\includegraphics[scale=0.55]{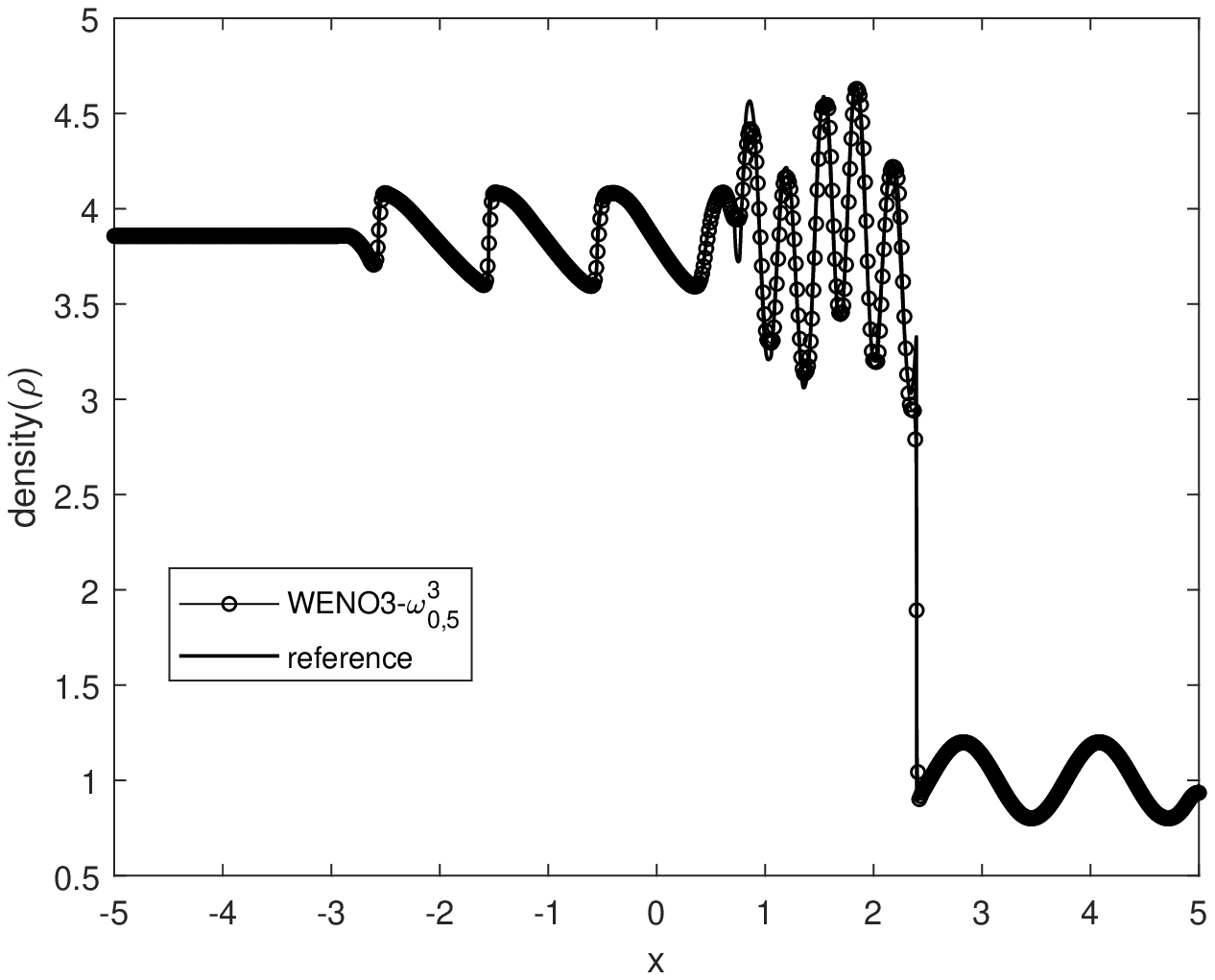} &
		\includegraphics[scale=0.55]{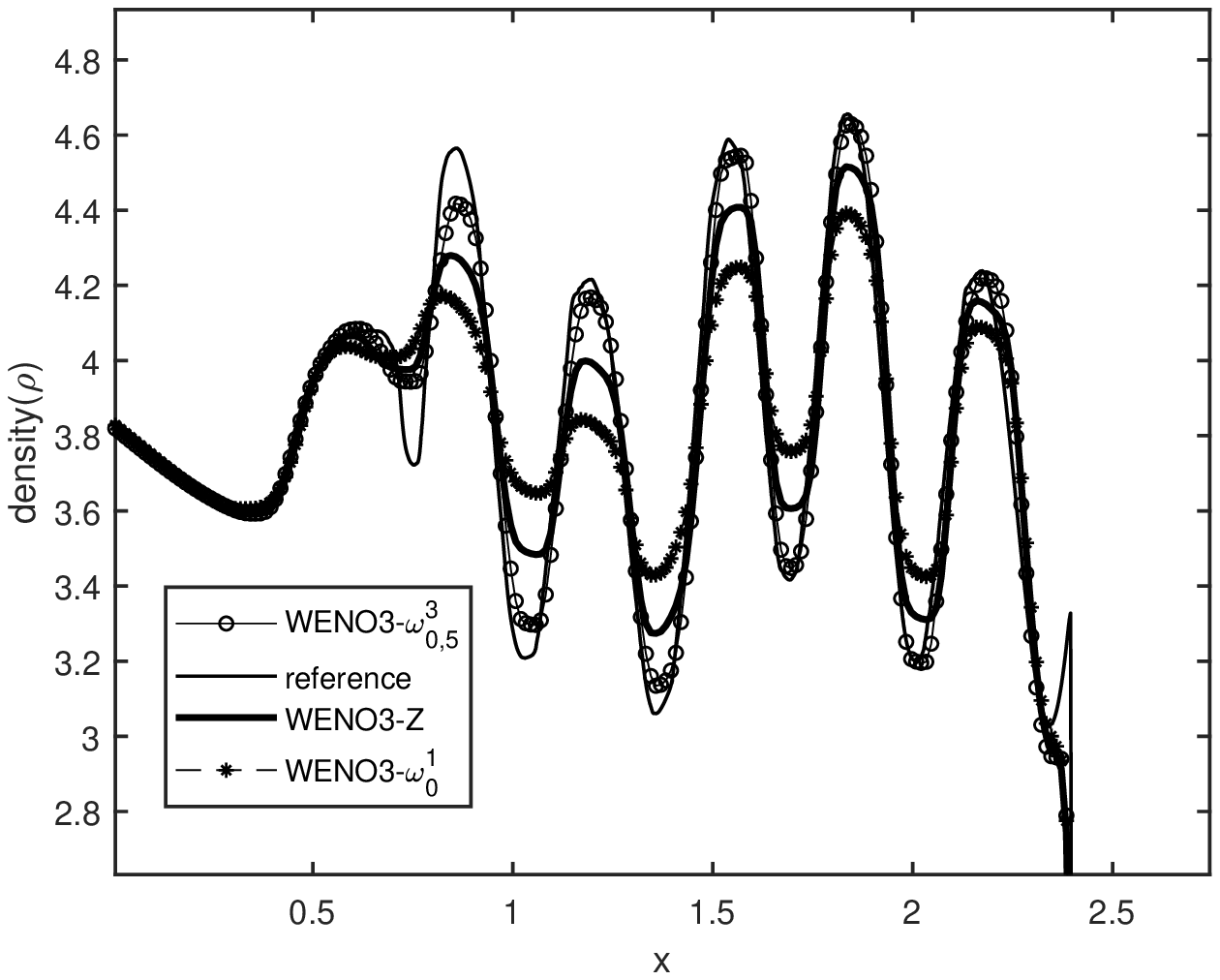}\\
		\textbf{a}&\textbf{b}
	\end{tabular}
	\caption{Solution of Shu-Osher test \eqref{shu-osher} using weight \textbf{(a)} WENO3-$\omega_{0,5}^3$ \textbf{(b)} comparison of WENO3-$\omega_{0,5}^3$ with WENO-$YC$ and WENO3-$\omega_{0}^1$. Solution is given at $t=1.8$ in the domain $[-5,5]$ with $800$ grids and $CFL=0.25$.}
	\label{fig:fig-5}
\end{figure}

The Shu-Osher test \cite{toro2013riemann} is the one-dimensional Euler system of equation \eqref{euler} with initial condition 
\begin{equation}\label{shu-osher}
(\rho,~u,~p)(x,0)=
\begin{cases}
(3.857143,~2.629369,~10.33333)~~~~~-5\le x<-4,\\
(1+0.2\sin(5x),~0,~1)~~~~~~~~~~~~~~~~~-4\le x\le 5,
\end{cases}
\end{equation}
in the computational domain $[-5,5]$and is run upto time $t=1.8$. Solution obtained using weight $\omega_{0,5}^3$ is given in In Figure \ref{fig:fig-5}(a) and compared in \ref{fig:fig-5}(b). It can be seen from Figure \ref{fig:fig-5}(b) that WENO3-$\omega_{0,5}^3$ resolve the smooth waves significantly better than WENO3-$YC$ and capture the shock crisply with no oscillations.
 \subsection{2D Riemann gas dynamic problem}
 Consider the two dimensional Euler equations of motion for gas dynamics
 \begin{equation}\label{2D}
 \begin{aligned} 
 \mathbf{w}_t+\mathbf{f}(\mathbf{w})_x+\mathbf{g}(\mathbf{w})_y&=0, \text{in}~ \mathbb R\times\mathbb R\times(0,\infty)\\ 
 \mathbf{w}&=\mathbf{w}_0 (x,y), ~\text{on}~ \mathbb R\times\mathbb R\times(t=0)
 \end{aligned}
 \end{equation}
 where $\mathbf{f}$ and $\mathbf{g}$ are fluxes in the $x$ and $y$ direction and $\mathbf{w}: \mathbb R\times \mathbb R\times (0,\infty)\to \mathbb R^m$ is the unknown with
 \begin{equation}
 \mathbf{w}=\left(\begin{array}{c}
 \rho\\
 \rho u\\
 \rho v\\
 E  \end{array}\right),~~\mathbf{f(w)}=\left(\begin{array}{c}
 \rho u\\
 \rho u^2+p\\
 \rho uv\\
 u(E+p)  \end{array}\right)~\text{and}~ \mathbf{g(w)}=\left(\begin{array}{c}
 \rho v\\
 \rho uv\\
 \rho v^2+p\\
 v(E+p) \end{array}\right)
 \end{equation}
 Total energy $E$ an the pressure $p$ are related by the follwing equation
 \begin{equation}
 p=(\gamma -1)(E-\frac{1}{2} \rho(u^2+v^2))
 \end{equation}
 Here $\rho,u,v,p$ and $E$ are density, components of velocity in the $x$ and $y$ coordinate directions,pressure, and total energy respectively.$\mathbf{w}$ is the vector of conservative variables,$\mathbf{f(w)}$ and $\mathbf{g(w)}$ is $x$ and $y$  direction wise flux component respectively.

\subsubsection{ Example 1:} \label{ex2}
\begin{figure}[ht!]
\begin{tabular}{cc}
	\centering
	\includegraphics[scale=0.5]{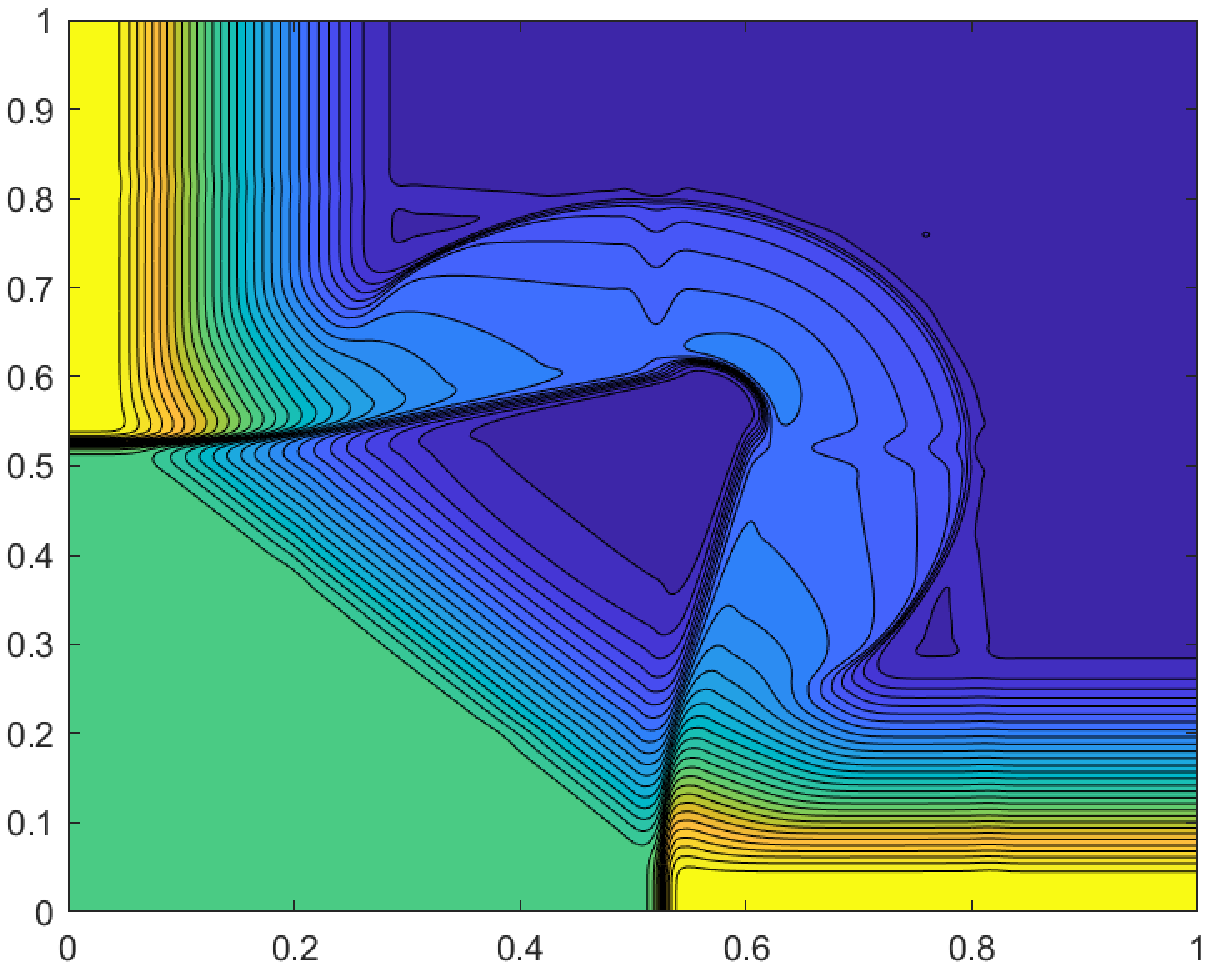} &
	\includegraphics[scale=0.5]{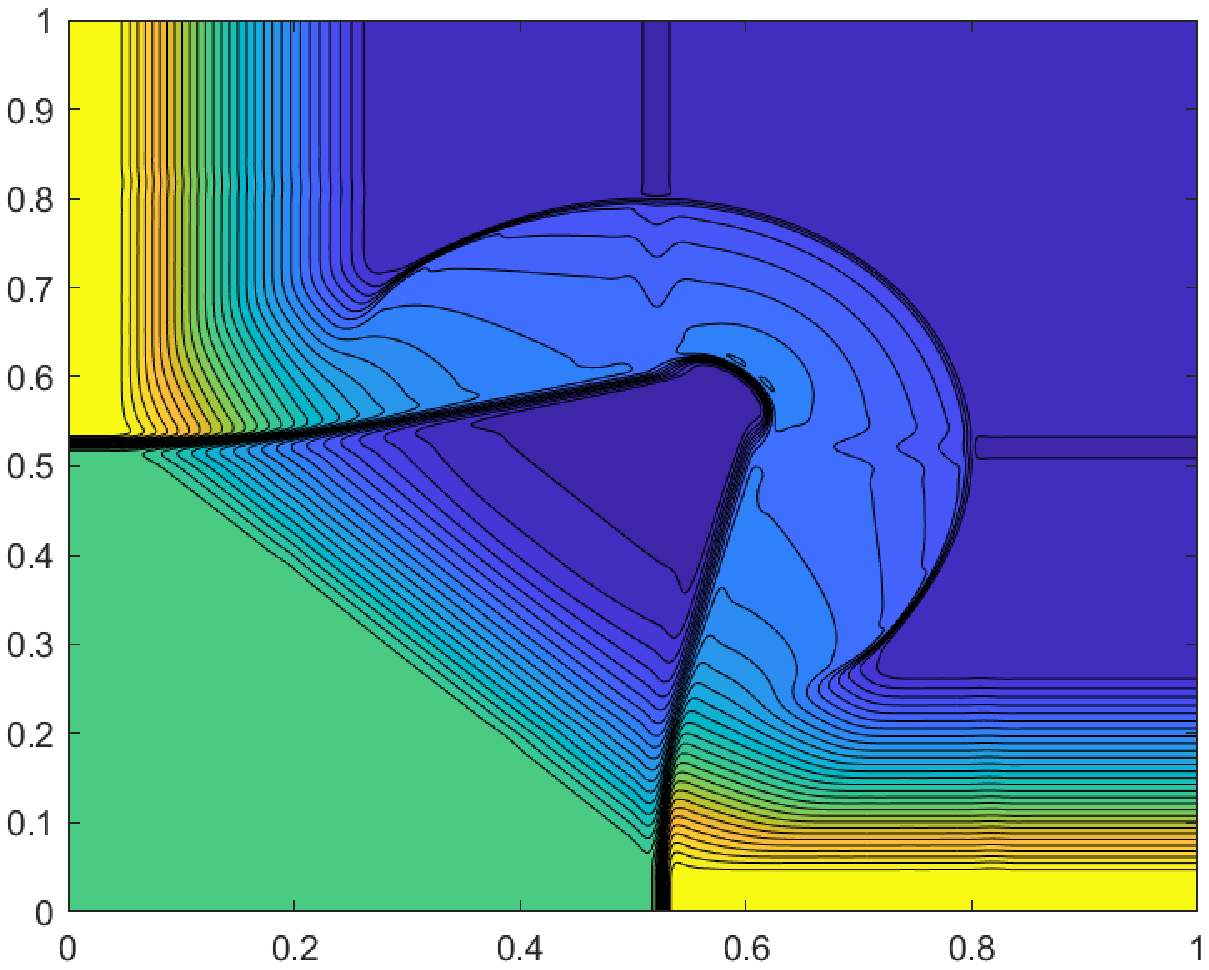}\\
	\textbf{a}&\textbf{b}
\end{tabular}
	\caption{Density profile of $2D$ example \ref{ex2} with $30$ contour lines by\textbf{(a)} WENO3-$\omega_{0}^1$ and \textbf{(b)}  WENO3-$\omega_{0,5}^3$ at $t=0.25$ using $CFL=0.45$ with $400\times 400$ grid points. }
	\label{fig:fig-8d}
\end{figure}

	The two-dimensional Riemann problem of gas dynamic \cite{kurganov2002solution} is defined by initial configuration as
	\begin{equation}
	(p,~\rho,~u,~v)(x,y,0)=
	\begin{cases}
	(0.4,~0.5197,~0.1,~0.1)~~~~~~\text{if}~x>0.5~\text{and}~y>0.5,\\
	(1,~1,-0.6259,~0.1)~~~~~~~~~\text{if}~x<0.5~\text{and}~y>0.5,\\
	(1,~0.8,~0.1,~0.1)~~~~~~~~~~~~\text{if}~x<0.5~\text{and}~y<0.5,\\
	(1,~1,~0.1,-0.6259)~~~~~~~~~\text{if}~x>0.5~\text{and}~y<0.5,
	\end{cases}
	\end{equation}
	The numerical solution is computed on the computational  square domain $[0,1]\times[0,1]$ with Dirichlet boundary conditions. The sqaure is divided into four quadrants by lines $x=0.5$ and $y=0.5$.And initial data as constant states are defined on each of the four quadrants and evolve upto time $t=0.25$ with $CFL=0.45$ for a grid $400\times 400$. Numerical result of WENO3-$\omega_{0}^1$ and WENO3-$\omega_{0,5}^3$ are given in Figure \ref{fig:fig-8d}.	

\subsubsection{Example 2:}\label{ex1}
\begin{figure}[htb!]
	\begin{tabular}{cc}
	\centering
		\includegraphics[scale=0.5]{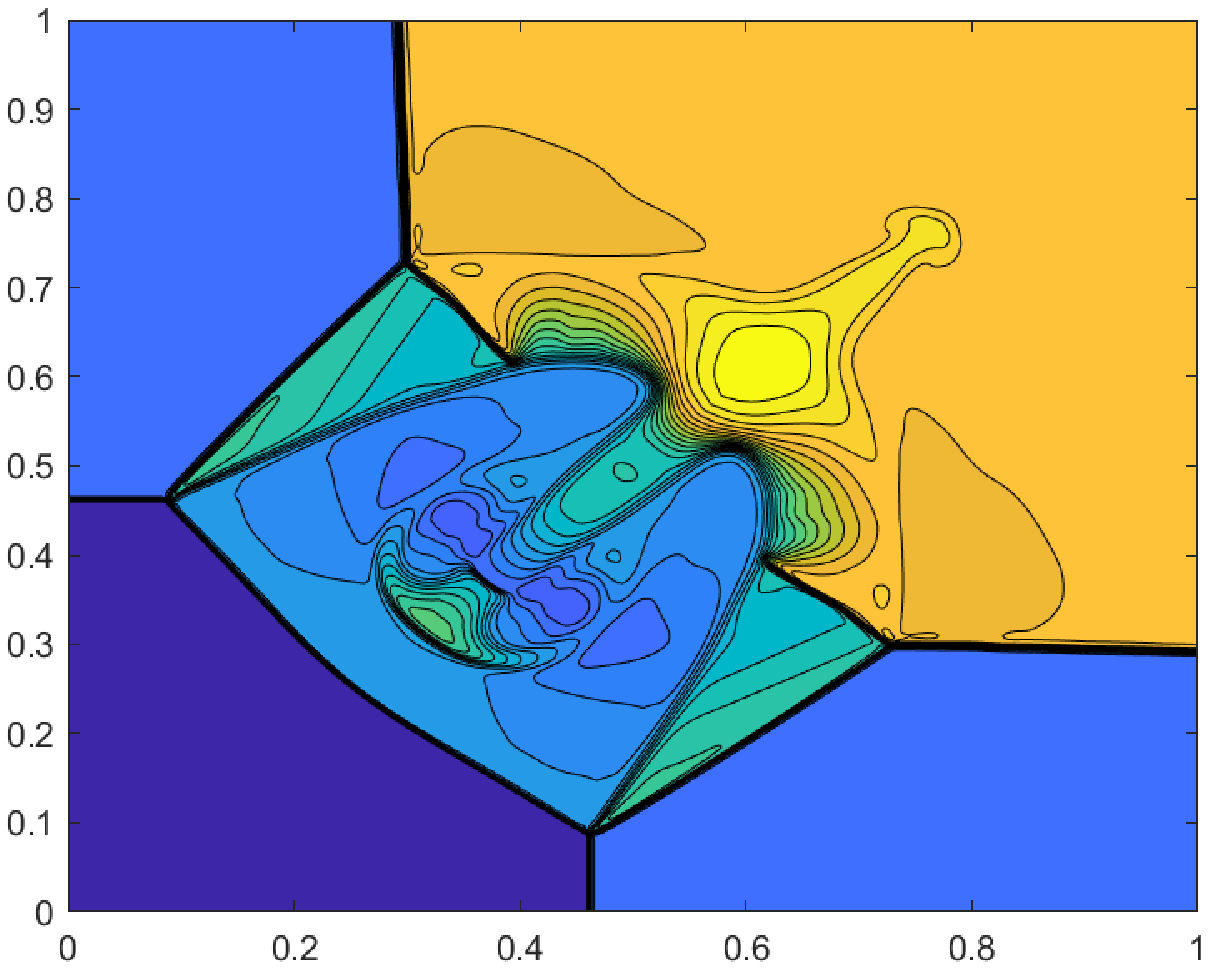} &
		\includegraphics[scale=0.5]{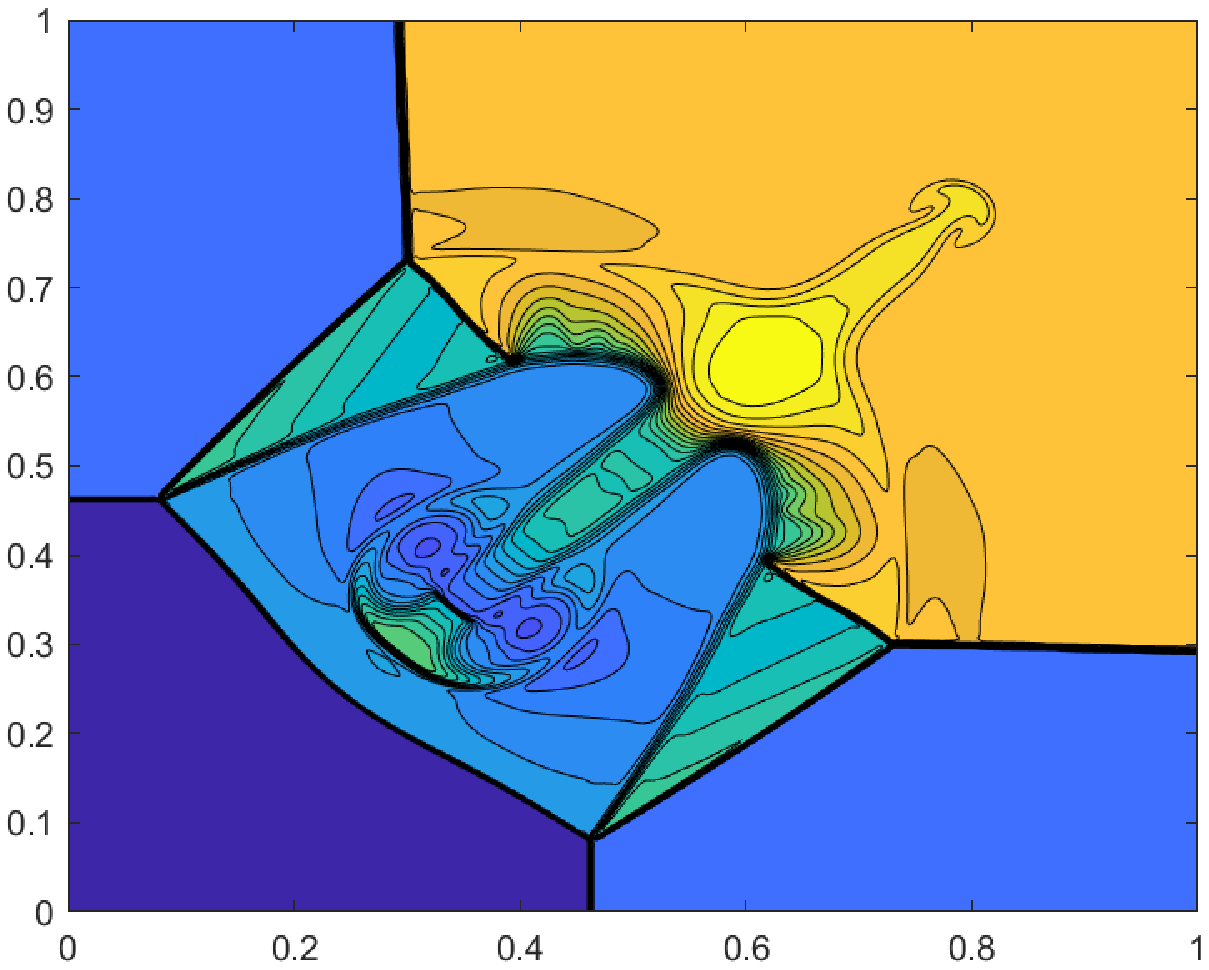}\\
		\textbf{a}&\textbf{b}
	\end{tabular}
	\caption{Density profile of $2D$ example \ref{ex1} with $30$ contour lines by \textbf{(a)} WENO3-$\omega_{0}^1$ and \textbf{(b)}  WENO3-$\omega_{0,5}^3$ at $t=0.8$ using $CFL=0.15$ with $400\times 400$ grid points. }
	\label{fig:fig-7d}
\end{figure}

This two-dimensional Riemann problem from \cite{schulz1993numerical} is defined by initial constant which is divided by the lines $x=0.8$ and $y=0.8$ as
\begin{equation}
(p,~\rho,~u,~v)(x,y,0)=
\begin{cases}
(1.5,~1.5,~0,~0)~~~~~~~~~~~~~~~~~~~\text{if}~0.8\le x\le 1~,~0.8\le y\le 1,\\
(0.3,~0.5323,~1.206,~0)~~~~~~~~~~\text{if}~0\le x <0.8~,~0.8\le y\le 1,\\
(0.029,~0.138,~1.206,~1.206)~~~\text{if}~0\le x<0.8~,~0\le y<0.8,\\
(0.3,~0.5323,~0,~1.206)~~~~~~~~~\text{if}~0.8< x\le 1~,~0\le y\le 0.8,
\end{cases}
\end{equation}
The numerical solution is computed on the computational  square domain $[0,1]\times[0,1]$ with Dirichlet boundary conditions. The sqaure is divided into four quadrants by lines $x=0.8$ and $y=0.8$.And initial data as constant states are defined on each of the four quadrants and evolve upto time $t=0.8$ with $CFL=0.15$ for a grid $400\times 400$. Numerical results of WENO3-$\omega_{0}^1$ and WENO3-$\omega_{0,5}^3$ are given in Figure \ref{fig:fig-7d}.	
\subsubsection{Explosion problem}
\begin{figure}[htb!]
	\begin{tabular}{cc}
		\centering
		\includegraphics[scale=0.5]{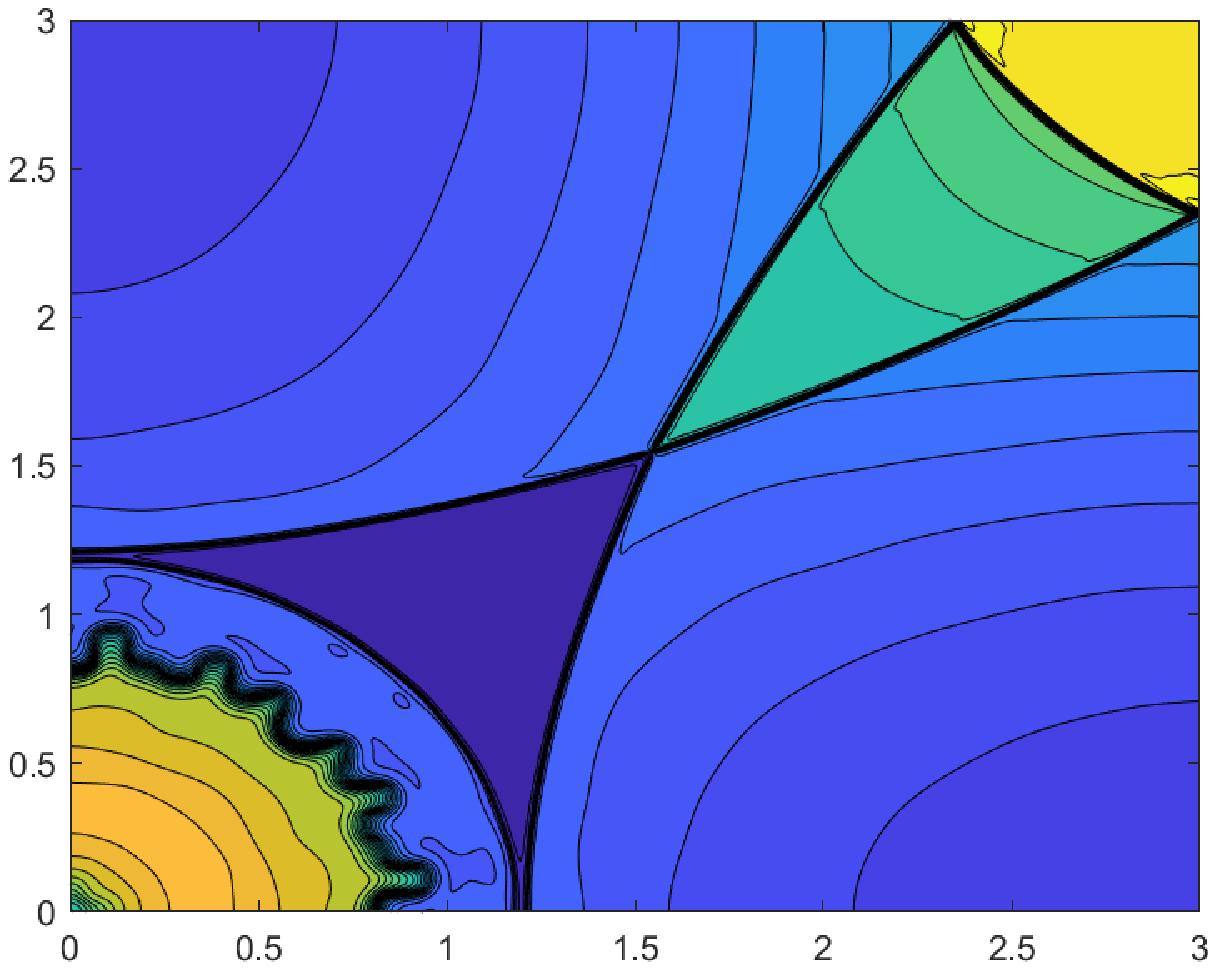} & \includegraphics[scale=0.5]{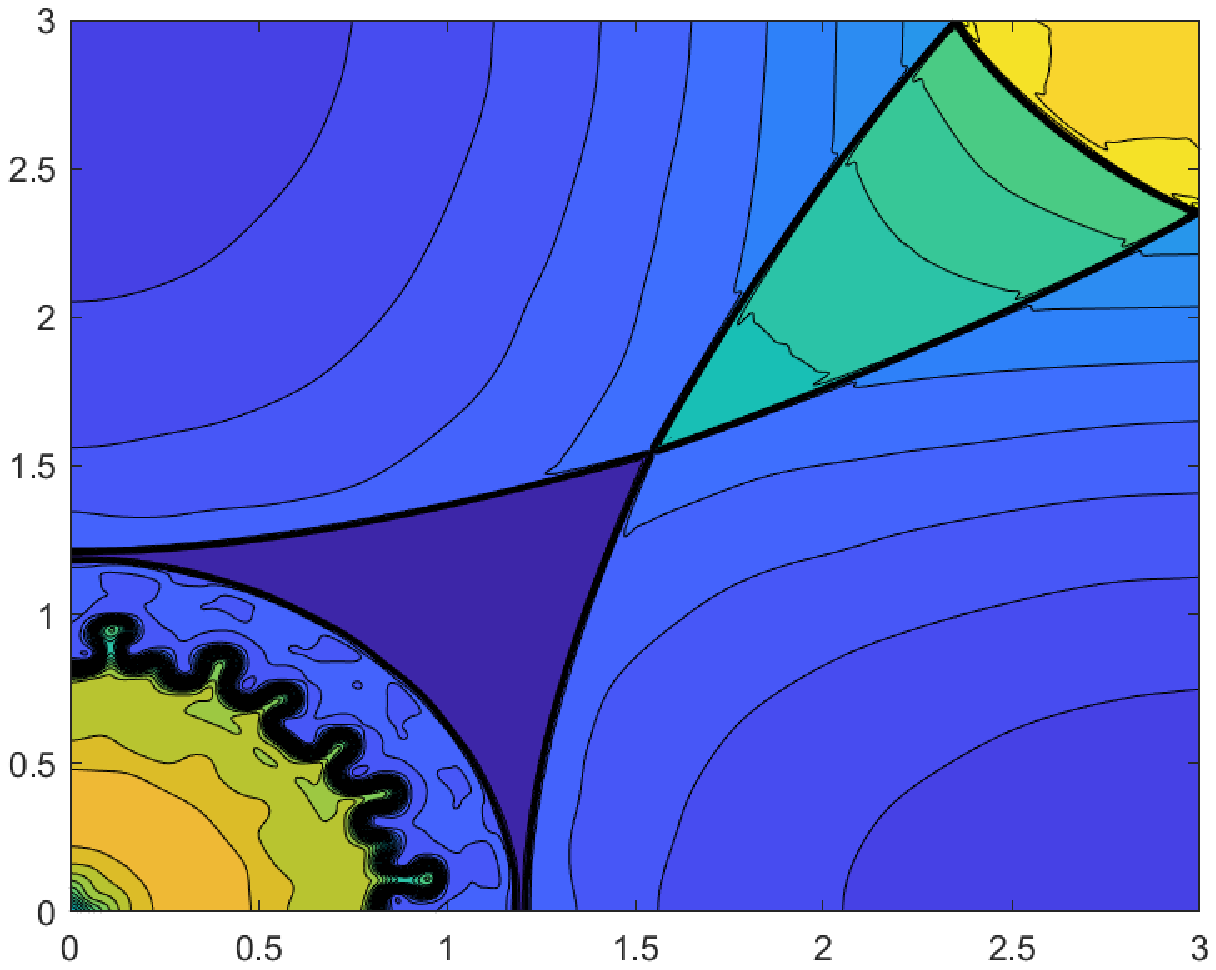} \\            \textbf{a}             &                        \textbf{b}
	\end{tabular}
	\caption{Results for Explosion problem  \textbf{(a)} WENO3-$\omega_{0}^1$, \textbf{(b)} WENO3-$\omega_{0,5}^3$ with $400\times 400$ grid to the time $t=3.2$ and $CFL=0.15$. }
	\label{fig:fig-10d}
\end{figure}

The explosion problem proposed in \cite{toro2013riemann} is a circularly symmetric $2D$ problem with initial circular region of higher density and higher pressure. In particular we set the center of the circle to the origin, its radius to $0.4$ and compute on a quadrant $(x,y)\in (0,1.5)\times(0,1.5)$. Inside the circle initial data are $p=1,~\rho=1,~u=0,~v=0$  and outside it is $p=0.1,~\rho=0.125,~u=0,~v=0$. i.e. the gas is initially at rest and its gas constant is 
$\gamma=\frac{1}{4}$. This problem (evolution of unstable contact at later times) is sensitive
to perturbations of the interface and as noted in \cite{toro2013riemann} for the cells which are crossed by the initial
interface circle one needs to use area weighted initial density and pressure. In Figure \ref{fig:fig-10d}, the numerical results for explosion problem are given.
 \subsubsection{Implosion problem}
\begin{figure}[htb!]
 \begin{tabular}{cc}
 	\centering
 	\includegraphics[scale=0.5]{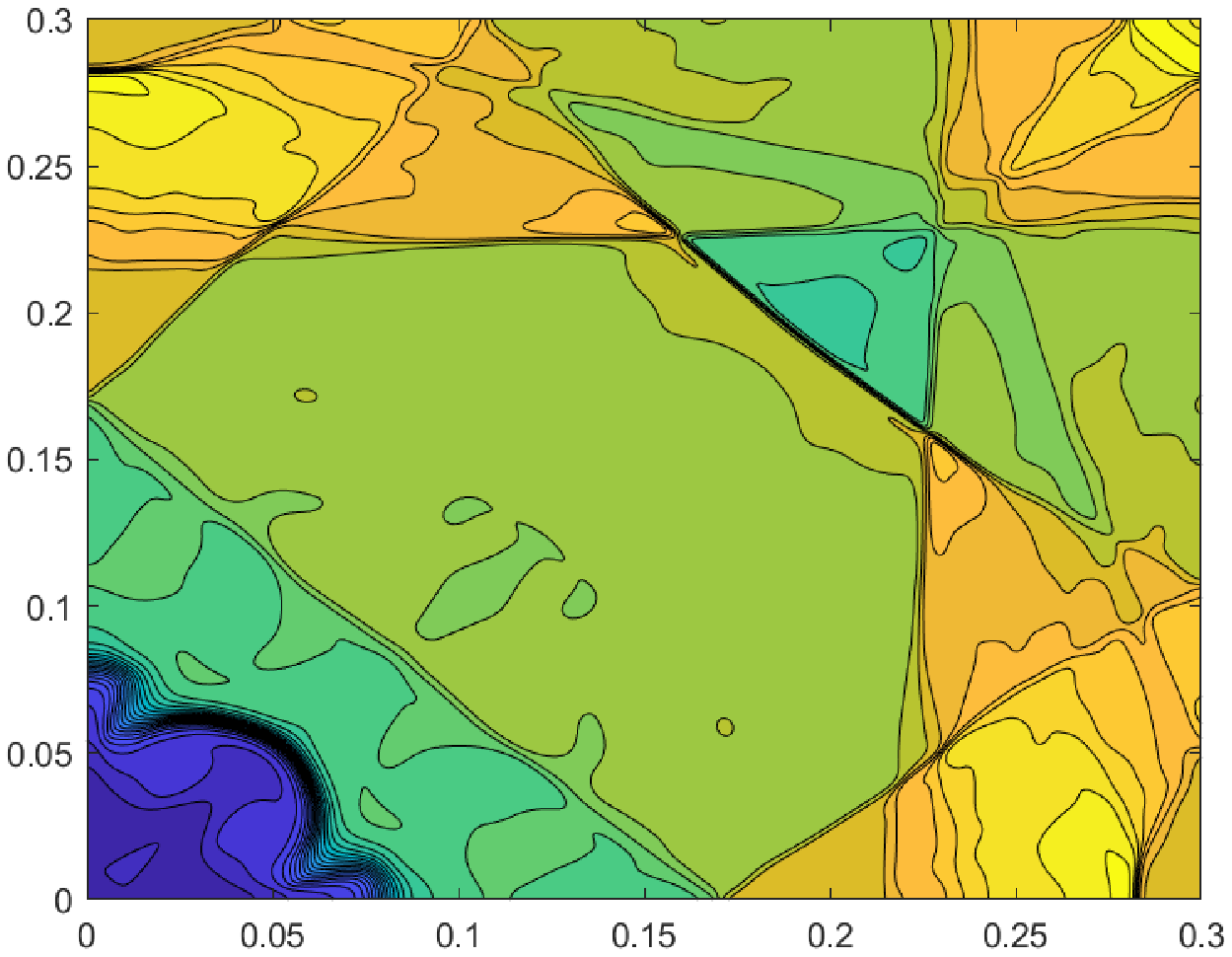} &
 	\includegraphics[scale=0.5]{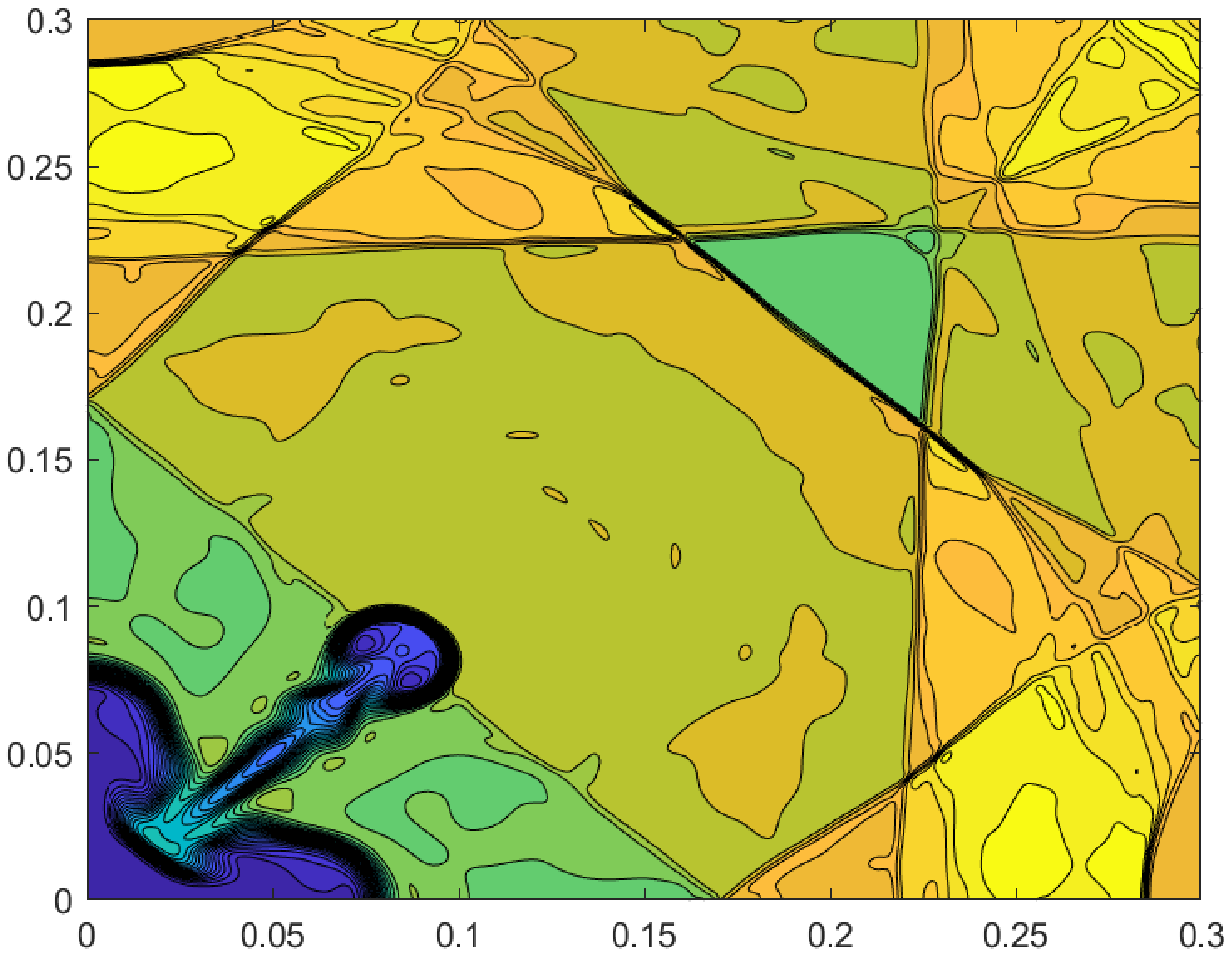}\\
 	\textbf{a}&\textbf{b}
 \end{tabular}
 	\caption{Results for Implosion problem \textbf{(a)} WENO3-$\omega_{0}^1$, \textbf{(b)} WENO3-$\omega_{0,5}^3$ with $400\times 400$ grid to the time $t=2.5$ and $CFL=0.25$. }
 	\label{fig:fig-9d}
 \end{figure}

 This shock problem has been presented in \cite{hui1999unified}.In this problem the gas is placed in a square box.The gas has initialy smaller density and pressure  inside a smaller square centered at the center of the box than the rest of the box.We use the box $(x,y)\in (-0.3,0.3)\times(-0.3,0.3)$ and the smaller square with corners at $(\pm 0.15,0),(0,\pm 0.15)$. we did the computation only in the  upper right quadrant $(x,y)\in (0,0.3)\times(0,0.3)$ of the box with diamond corner box $|x|+|y|<0.15$. Initial data inside the diamond corner box are $p=0.14,~\rho=0.125,~u=0,~v=0$ and outside are $p=1.0,~\rho=1,~u=0,~v=0$. i.e. initial velocities are zero. The gas constant is $\gamma=\frac{1}{4}$. On all four boundaries Reflecting boundary conditions are used. In Figure \ref{fig:fig-9d}, the numerical results for implosion problem are given.
 \\
It can be observed from above four 2D Reimann Problem examples that, both the schemes WENO3-$\omega_{0}^1$, and WENO3-$\omega_{0,5}^3$ resolves the flow feature however scheme WENO3-$\omega_{0,5}^3$ gives significantly better resolution to discontinuities compared to scheme WENO3-$\omega_{0}^1$.

 \pagebreak
\section{Conclusion \& Future Work}\label{sec5}
A new framework is given to construct non-linear weight by utilizing necessary conditions for non-oscillatory WENO3 reconstruction. A characterization non-linear weights for constructing third order WENO schemes is done and various nonlinear weights using weight limiting functions are given. Computational results are given and in some cases compared with WENO3 schemes using other established weights. These results show third order accuracy of the resulting schemes and ability to resolve smooth as well discontinuous region of the solution. 
The work to propose similar framework for fifth order WENO scheme is under progress and partial computational results are obtained. However, an extensive theoretical investigation is further required to get suitable necessary conditions. This work will be reported separately. \\\\
{\bf Acknowledgment:} Authors acknowledge the Science and Engineering Board, New Delhi, India for providing necessary financial support through funded projects File No. EMR/2016/ 000394 and MTR/2017/000187.

\bibliographystyle{unsrt}
\bibliography{ref}
\end{document}